\documentclass[11pt,reno]{amsart}  
\usepackage{amssymb}
\usepackage{amsmath}
\usepackage{enumitem}
\usepackage{mathrsfs}
\usepackage[english]{babel}
\usepackage[all]{xy}
\usepackage{hyperref}
\usepackage{marginnote}
\usepackage{stmaryrd}
\usepackage{fullpage}

\RequirePackage{mathrsfs}

\newtheorem{theorem}{Theorem}[section]
\newtheorem{lemma}[theorem]{Lemma}
\newtheorem{proposition}[theorem]{Proposition}
\newtheorem{corollary}[theorem]{Corollary}
\newtheorem{conjecture}[theorem]{Conjecture}

\theoremstyle{definition}
\newtheorem*{ack}{Acknowledgements}
\newtheorem*{con}{Conventions}
\newtheorem*{motivation}{Motivation}
\newtheorem{remark}[theorem]{Remark}

\newtheorem{definition}[theorem]{Definition}

\numberwithin{equation}{section} \numberwithin{figure}{section}

\DeclareMathOperator{\Spec}{Spec}
\DeclareMathOperator{\Cov}{Cov}

\DeclareMathOperator{\id}{id}

\newcommand{\cl}[1]{\mathrm{cl}\left(#1\right)}

\title{Products of varieties with many integral points}
\author{Cedric Luger}
\address{Cedric Luger\\
Institut f\"{u}r Mathematik\\
Johannes Gutenberg-Universit\"{a}t Mainz\\
Staudingerweg 9, 55099 Mainz\\
Germany.}
\email{celuger@uni-mainz.de}



\begin{document}

\begin{abstract}
    Corvaja and Zannier asked whether a smooth projective integral variety with a dense set of rational points over a number field satisfies the weak Hilbert property.
    We introduce an extension of the weak Hilbert property for schemes over arithmetic base rings by considering near-integral points, extending Corvaja--Zannier's question beyond the projective case.
    Building on work of Bary-Soroker--Fehm--Petersen and Corvaja--Demeio--Javanpeykar--Lombardo--Zannier, we prove several properties of this more general notion, in particular its persistence under products.
    We also answer positively Corvaja--Zannier's question for all algebraic groups over finitely generated fields of characteristic zero.
\end{abstract}

\maketitle

\thispagestyle{empty}

\section{Introduction}

An integral variety $X$ over a field $k$ of characteristic $0$ has the \emph{Hilbert property (over $k$)} if $X(k)$ is not thin in $X$, as defined by Serre \cite[§3]{SerreTopicsGalois}.
Applying the Chevalley--Weil Theorem, Corvaja--Zannier \cite[Theorem~1.4]{CZHP} showed that a smooth projective integral variety with the Hilbert property over a number field is algebraically simply connected.
In light of this, it is more natural to study the notion of the \emph{weak Hilbert property}, in which one restricts to ramified covers in the definition of thin sets, see \cite[Definition~1.2]{CDJLZ}.
(By a \emph{ramified cover}, we mean a finite surjective morphism that is not unramified.)

\begin{definition}
    Let $X$ be a normal integral variety over a field $k$ of characteristic $0$. A subset $\Sigma \subseteq X(k)$ is \emph{strongly thin (in $X$ or in $X(k)$)} if there exists
    a non-negative integer $n$ and, for every $i=1,\ldots,n$, a normal integral variety $Y_i$ over $k$ and a ramified cover $\pi_i \colon Y_i \to X$ such that
    $\Sigma \setminus \bigcup_{i=1}^n \pi_i(Y_i(k))$ is not dense in $X$.
\end{definition}

We restrict ourselves in the above definition to normal varieties to ensure that ramification is well-behaved (in particular, a cover is either ramified or étale, see \cite[Lemma~2.3]{CDJLZ}).
Note that, if $Y_i$ is not geometrically integral, then $Y_i(k) = \emptyset$, so we may always assume $Y_i$ to be geometrically integral.
We are interested in studying strongly thin subsets of quasi-projective varieties given by sets of integral points. To do this, we use the language of arithmetic schemes.

If $S$ is an irreducible scheme with function field $K(S)$, we say that an $S$-scheme is an \emph{arithmetic scheme over $S$} if it is dominant separated of finite type over $S$ and the generic fibre $\mathcal{X}_{K(S)}$ is normal and integral.
By an \emph{arithmetic scheme} we mean an arithmetic scheme over $\mathbb{Z}$.

In order to extend Corvaja--Zannier's definition of the weak Hilbert property to arithmetic schemes, we consider near-integral points, as defined by Vojta \cite[Definition~4.3]{VojtaLangExc} (see also \cite[Definition~3.11]{JavanpeykarLVConjectures}).

\begin{definition}
    Let $S$ be an integral noetherian scheme and let $\mathcal{X} \to S$ be a finite type morphism.
    We say that a $K(S)$-point $P$ of $\mathcal{X}$ is \emph{near-integral} (or \emph{near-$S$-integral}) if there is a dense open $U \subseteq S$ whose complement is of codimension at least two such that $P$ extends to a morphism $U \to \mathcal{X}$.
    If $T \to S$ is a morphism of integral noetherian schemes, then we define $\mathcal{X}(T)^{(1)} \subseteq \mathcal{X}_T(K(T)) = \mathcal{X}(K(T))$ to be the set of near-$T$-integral points of $\mathcal{X}_T \to T$.
\end{definition}

Clearly, a morphism $\mathcal{X} \to \mathcal{Y}$ of schemes of finite type over $S$ induces a map $\mathcal{X}(S)^{(1)} \to \mathcal{Y}(S)^{(1)}$.
If $S = \Spec R$, we also write $\mathcal{X}(R)^{(1)} = \mathcal{X}(S)^{(1)}$ and refer to near-$S$-integral points as near-$R$-integral points.
For example, if $\mathcal{X}$ is a proper arithmetic scheme over a normal $\mathbb{Z}$-finitely generated integral domain $R$ with fraction field $k$ of characteristic zero, then the set of near-$R$-integral points of $\mathcal{X}$ equals its set of $k$-points.
If $R$ is $1$-dimensional, e.g. $R=\mathcal{O}_{K,S}$ is the ring of $S$-integers of a number field $K$ for some finite set of finite places $S$ of $K$, then the near-$R$-integral points coincide with the $R$-integral points.

The weak Hilbert property can be naturally extended to arithmetic schemes as follows; such a notion was also studied in \cite{CocciaCubic} over $S$-integers of number fields.

\begin{definition}\label{Def:WHP}
    Let $S$ be a normal noetherian integral scheme with function field $K(S)=k$ of characteristic $0$ and let $\mathcal{X}$ be an arithmetic scheme over $S$.
    We say that $\mathcal{X}$ has the \emph{weak Hilbert property (over $S$)} if $\mathcal{X}(S)^{(1)}$ is not strongly thin in $\mathcal{X}_k$. 
\end{definition}

We will primarily study the weak Hilbert property over normal integral arithmetic base schemes; for Hilbert irreducibility type results over arbitrary rings or fields see for example \cite{BSFP14, BSFP23BaseChange, BSFP23Abelian, FriedJarden, LangFundamentals}.
A similar extension of the (usual) Hilbert property to arithmetic schemes is given in \cite{LugerArithmeticHP}.

Assuming the base scheme is arithmetic, the Hilbert property and weak Hilbert property differ only in the presence of finite étale covers.
In fact, Corvaja and Zannier show in \cite[Theorem~1.4]{CZHP} that a smooth projective variety with the Hilbert property over a number field is algebraically simply connected (i.e., admits no non-trivial finite étale cover); their result also holds for normal arithmetic schemes \cite[Theorem~1.8]{LugerArithmeticHP}.

Our goal is to show that the weak Hilbert property is closed under products (see Theorem~\ref{Thm:FibrationProducts} for a more general result pertaining to fibrations):

\begin{theorem}\label{Thm:ProductsWHP}
    Let $S$ be a regular integral arithmetic scheme
    and let 
    $\mathcal{X}$ and $\mathcal{Y}$ be quasi-projective arithmetic schemes over $S$.
    If both $\mathcal{X}$ and $\mathcal{Y}$ satisfy the weak Hilbert property over $S$, then $\mathcal{X} \times_S \mathcal{Y}$ also satisfies the weak Hilbert property over $S$.
\end{theorem}

Note that in Theorem~\ref{Thm:ProductsWHP}, the product $\mathcal{X} \times_S \mathcal{Y}$ is again an arithmetic scheme over $S$.
Indeed, since the generic fibres $\mathcal{X}_{K(S)}$ and $\mathcal{Y}_{K(S)}$ have a dense set of $K(S)$-rational points by the weak Hilbert property, they are geometrically irreducible by \cite[Tag~0G69]{stacks-project} and thus normal geometrically integral varieties (being geometrically reduced is automatic by the characteristic of $K(S)$ being $0$ since $S$ is dominant of finite type over $\Spec \mathbb{Z}$). Therefore, their product is also normal and (geometrically) integral over $K(S)$.
If $\mathcal{X}$ is an arithmetic scheme over $S$, then $\mathcal{X}$ has exactly one irreducible component $\mathcal{X}_0$ that dominates $S$ (as the generic fibre is irreducible). Moreover, by the irreducibility of $S$, we have $\mathcal{X}(S)^{(1)} = \mathcal{X}_0(S)^{(1)}$, so when studying the weak Hilbert property, we are usually only interested in irreducible arithmetic schemes.
However, if $\mathcal{X}$ and $\mathcal{Y}$ are irreducible arithmetic schemes with geometrically irreducible generic fibre over an irreducible arithmetic scheme $S$, then $\mathcal{X} \times_S \mathcal{Y}$ is an arithmetic scheme over $S$, though not necessarily irreducible. For this reason, we don't require arithmetic schemes to be irreducible in the definition.
Note that this issue does not arise if $\mathcal{X}$ and $\mathcal{Y}$ are flat over $S$, e.g., if they are irreducible and $\dim S = 1$.

We also prove some further basic properties. Namely,
the weak Hilbert property
is a birational invariant (Proposition~\ref{Prop:WHPBirational}), it
is inherited by smooth proper images (Proposition~\ref{Prop:WHP_descends_smooth_proper}),
and finite étale covers inherit the weak Hilbert property, up to possibly extending the base (Proposition~\ref{Prop:PersistencePotentialWHP}).

Theorem~\ref{Thm:ProductsWHP} generalizes a similar product theorem for smooth proper varieties over a finitely generated field $k$ of characteristic $0$
(see \cite[Theorem~1.9]{CDJLZ}).
In fact, our proof of Theorem~\ref{Thm:ProductsWHP} follows the proof of \emph{loc.\ cit.}\ closely, with the main novelty being the elimination of the properness assumption in several places.
The lack of properness complicates especially the structure of vertically ramified covers, see Lemma~\ref{Lem:VerticalRamification} for the key technical generalization of \cite[Lemma~2.17]{CDJLZ}.

Over arbitrary fields, the (usual) Hilbert property was shown to satisfy a similar product property by Bary-Soroker--Fehm--Petersen \cite[Corollary~3.4]{BSFP14}
(see \cite[Theorem~1.5]{LugerArithmeticHP} for a generalization to arithmetic schemes),
answering positively a question of Serre.
(A new proof is given in \cite[Corollary~3.5]{BSFP23BaseChange}, and over
number fields this also follows from \cite[Lemma~8.12]{HarpazWittenberg}.)
An additional difficulty in proving the product theorem for the weak Hilbert property is the possible presence of a vertically ramified cover of a product, which is where the arithmeticity of the base scheme is required.

We were first led to investigate the weak Hilbert property by a conjecture of Corvaja--Zannier  \cite[Question-Conjecture~2]{CZHP} which predicts that a smooth projective variety with a dense set of rational points over a number field should have the weak Hilbert property, after a suitable extension of the base field.
A similar statement should hold for $S$-integral points of number fields (see for example \cite[Question-Conjecture~1bis]{CZHP}).
The following conjecture is a natural extension of the conjectures of Corvaja--Zannier to the setting of quasi-projective arithmetic schemes:

\begin{conjecture}\label{Conj:CZ_arithmetic}
    Let $R$ be a normal $\mathbb{Z}$-finitely generated integral domain with fraction field $k$ of characteristic zero and let $\mathcal{X}$
    be a regular quasi-projective arithmetic scheme over $R$. If $\mathcal{X}(R)^{(1)}$ is dense in $\mathcal{X}_k$, then there is a
    dominant generically finite morphism $\Spec S \to \Spec R$ of normal noetherian integral schemes such that $\mathcal{X}_S$ has the weak Hilbert property over $S$.
\end{conjecture}

Note that, as before, the density of $\mathcal{X}(R)^{(1)}$ implies that $\mathcal{X}_k$ is geometrically integral. Therefore, $\mathcal{X}_{K(\Spec S)}$ is a normal integral variety, showing that $\mathcal{X}_S$ is an arithmetic scheme over $\Spec S$.

It follows from Zannier's results \cite[Theorem~1 and Theorem~2]{ZannierDuke}
that powers of elliptic curves and group schemes of the form $\mathbb{G}_{\mathrm{a}}^n \times \mathbb{G}_{\mathrm{m}}$
over a ring of $S$-integers of a number field satisfy Conjecture~\ref{Conj:CZ_arithmetic}
if they have a non-degenerate $k$-point (a point that generates a dense cyclic subgroup).
We prove a generalization of this statement for group schemes.
More explicitly, we show the following:
\begin{theorem}\label{Thm:AlgGroupsSatisfyConjCZ}
    Let $R$ be a normal $\mathbb{Z}$-finitely generated integral domain with fraction field $k$ of characteristic zero and
    let $\mathcal{G}$ be a connected finite type group scheme over~$R$.
    If $\mathcal{G}(k)$ is dense in $\mathcal{G}$,
    then there exists a dense open $\Spec R' \subseteq \Spec R$
    such that $\mathcal{G}(R')$ is not strongly thin in $\mathcal{G}_k$.
    In particular, $\mathcal{G}_{R'}$ has the weak Hilbert property over $R'$.
\end{theorem}
Note that Theorem~\ref{Thm:AlgGroupsSatisfyConjCZ} is ``optimal'' in the sense that we can not expect $R'$ to be equal to $R$;
consider for example $R=\mathbb{Z}$ and $\mathcal{G} = \mathbb{G}_{\mathrm{m}, \mathbb{Z}}$.

The proof of Theorem~\ref{Thm:AlgGroupsSatisfyConjCZ} relies on work of Liu \cite{FeiLiu} and Corvaja \cite{CorvajaFixedPoints} on linear algebraic groups and Corvaja--Demeio--Javanpeykar--Lombardo--Zannier \cite{CDJLZ} on abelian varieties.
It follows from \emph{loc.\ cit.}\ that Theorem~\ref{Thm:AlgGroupsSatisfyConjCZ} holds for $R'=R$ under the additional assumption that $\mathcal{G}(R)$ is dense and either $\mathcal{G}_k$ is linear or $\mathcal{G}$ is an abelian group scheme (see also Proposition~\ref{Prop:FiniteIndexCosetManyRamCovers} and Proposition~\ref{Prop:FinGenDenseSubgroup}),
though it remains unknown at this point if this is also true for general $\mathcal{G}$.
In a nutshell, to prove Theorem~\ref{Thm:AlgGroupsSatisfyConjCZ}, we have to deduce non-liftability of rational points along finite collections of ramified covers from properties of (PB)-covers; see Section~\ref{Sec:AlgGroups} for details.

As a direct consequence, we obtain that every connected algebraic group with a dense set of rational points over a finitely generated field $k$ of characteristic $0$ satisfies the weak Hilbert property
(see Corollary~\ref{Cor:AlgGroupsConjCZ}).
Other instances of Conjecture~\ref{Conj:CZ_arithmetic} are proved in
\cite{CocciaCubic, CocciaDelPezzo, DemeioStreeter, DemeioStreeterWinter, Demeio1, Demeio2, GvirtzChenHuang, GvirtzChenMezzedimi, GvirtzChenMezzedimi2, JavanpeykarProductOfCurves, JavanpeykarNefTangentBundle, LoughranSalgado, NakaharaStreeter, Streeter}.

\begin{motivation}
The motivation for this paper is to make it possible to investigate sets of near-integral points on products of arithmetic schemes that are not necessarily given by a product of sets of near-integral points.
For example, if $\mathcal{X}$ and $\mathcal{Y}$ are arithmetic schemes over a suitable ring $R$ and $\mathcal{Z} \subseteq \mathcal{X} \times_R \mathcal{Y}$ is a closed subscheme, then Theorem~\ref{Thm:FibrationProducts} applies to near-$R$-integral points on the ``punctured'' scheme $(\mathcal{X} \times_R \mathcal{Y})\setminus \mathcal{Z}$.
In fact, the results of this paper are crucial in our recent work \cite{LugerPuncturedWHP} on the conjectures of Campana and Corvaja--Zannier for ``punctured'' algebraic groups.
\end{motivation}
 
\begin{ack}
We gratefully acknowledge Ariyan Javanpeykar for constant support and many helpful discussions.
We also thank Finn Bartsch, Manuel Blickle, and Georg Tamme for helpful discussions.
We gratefully acknowledge the referee for a careful reading of the paper and for numerous insightful comments, including the identification of several inaccuracies and valuable suggestions on how to address them.
This research was funded by the Studienstiftung des Deutschen Volkes (German Academic Scholarship Foundation).
We gratefully acknowledge support by the Deutsche Forschungsgemeinschaft (DFG, German Research Foundation) -- Project-ID 444845124 -- TRR 326.
Parts of the results contained in this work were obtained during the author's stays at the Radboud University Nijmegen. We are grateful for accomodation and an excellent working environment.
\end{ack}

\begin{con}
    A \emph{variety} over a field $k$ is a reduced separated scheme of finite type over $k$.
    A \emph{cover} is a finite surjective morphism of integral varieties, and
    a cover of normal integral varieties is
	\emph{ramified} if it is not unramified.
 	If $S$ is a scheme and $s\in S$ a point, we let $\kappa(s)$ denote the residue field of $s$.
 	If $S$ is irreducible with generic point $\eta$, we write $K(S) = \kappa(\eta)$ for the function field of $S$.
\end{con}

\section{Vertically ramified covers}

In \cite[Lemma~2.17]{CDJLZ} it is proved that, for smooth proper varieties $X, Y$, a cover $Z \to X\times Y$ ``vertically ramified'' over $X$ (that is, $Z \to X \times Y$ is unramified over $U\times Y$ for some dense open $U \subseteq X$) splits as a product of a ramified cover of $X$ and an unramified cover of $Y$, up to a finite étale cover. We provide a generalization that eliminates the properness assumptions.
For the following proof to work without assuming properness, we require $X(k)$ to be dense (instead of $U(k)\neq \emptyset$, as in \cite{CDJLZ}).

\begin{lemma}\label{Lem:VerticalRamification}
    Let $k$ be a field of characteristic zero, $X, Y, Z$ normal integral varieties over $k$, and $\pi \colon Z\to X \times Y$ a ramified cover. Assume that there exists a dense open subscheme $U \subseteq X$ such that $\pi$ is unramified over $U\times Y$ (i.e., $\pi$ is vertically ramified over $X$).
    Let $Y \to \overline{Y}$ be an open immersion into a normal integral proper variety $\overline{Y}$ over $k$, and let $\overline{Z} \to X \times \overline{Y}$ be the normalization of $X \times \overline{Y}$ in $Z$. Assume that the fibres of $\overline{Z} \to X$ are geometrically connected and $Z(k)$ is dense. Then there exists a commutative diagram
    \[
     \xymatrix{
        & \overline{Z'} \ar[dl] \ar[d] & Z' \cong X' \times Y' \ar[l] \ar[d] \ar[r] & Z \ar[d] \\
        X' \ar[dr] & X \times \overline{Y'} \ar[d] & X \times Y' \ar[l] \ar[d] \ar[r] & X \times Y \ar[d] \\
        & X \ar@{=}[r] & X \ar@{=}[r] & X
     }
    \]
    such that
    \begin{enumerate}
        \item $X', Y', Z', \overline{Y'}, \overline{Z'}$ are normal integral varieties over $k$;
        \item \label{Item:VerticalRam:Stein_is_ramified} $X' \to X$ is a ramified cover;
        \item $Y' \to Y$ and $Z' \to Z$ are finite étale;
        \item $Z'$ is a connected component of the fibre product $Z \times_Y Y'$;
        \item $\overline{Y'}$ is the normalization of $\overline{Y}$ in $Y'$;
        \item $\overline{Z'}$ is the normalization of $X \times \overline{Y'}$ in $Z'$;
        \item $\overline{Z'} \to X' \to X$ is the Stein factorization of $\overline{Z'} \to X$.
    \end{enumerate}
\end{lemma}
\begin{proof}
    Let $\overline{k}$ be the algebraic closure of $k$.
    Since $Z(k)$ is dense in $Z$ and thus $X(k)$ and $Y(k)$ are dense in $X$ and $Y$, respectively, it follows from \cite[Tag~0G69]{stacks-project} that $X_{\overline{k}}$ and $Y_{\overline{k}}$ are irreducible.
    Since $\overline{k}$ is algebraically closed of characteristic $0$,
    the étale fundamental group of $Y_{\overline{k}}$ is topologically finitely generated
    (see \cite[Théorème~2.3.1]{SGA7I_Raynaud} and \cite{HironakaResolution}),
    so there are only finitely many $Y_{\overline{k}}$-isomorphism classes of finite étale morphisms $V \to Y_{\overline{k}}$ of bounded degree.
    Since $\overline{Z}$ is normal and $k$ is of characteristic zero, we may shrink $U$ if necessary to assume that $(\overline{Z})_x$ is normal for every $x\in U(k)$. In particular, since $(\overline{Z})_x$ is geometrically connected, this implies that $Z_x$ is normal geometrically connected (and thus geometrically integral) for every $x\in U(k)$.
    Let $Z^1 \to Y_{\overline{k}}, \ldots, Z^n \to Y_{\overline{k}}$ be representatives of the classes of finite étale morphisms of degree $\deg(Z \to X\times Y)$.
    For $x\in X(k)$, let $\overline{x}\colon \Spec \overline{k} \to X$ be the induced geometric point.
    Define $\Sigma_i \subseteq U(k)$ to be the subset of $x\in U(k)$ such that (the finite étale morphism) $Z_{\overline{x}} = (Z_x)_{\overline{k}} \to (\{x\} \times Y)_{\overline{k}} \cong Y_{\overline{k}}$ is $Y_{\overline{k}}$-isomorphic to $Z^i$. Note that $U(k) = \Sigma_1 \cup \ldots \cup \Sigma_n$.

    Since $U(k)$ is dense in $X$, there exists an $i$ such that $\Sigma_i$ is dense in $X$;
    we define $\Sigma_i =: \Sigma$.
    Fix some point $x' \in \Sigma$ and define $Y' = Z_{x'}$.
    Then $Y'$ is a normal geometrically connected $k$-variety, $Y' \to Y$ is finite étale,
    and $\Sigma$ is a dense subset of $U(k)$ such that, for every $x\in \Sigma$, there exists a $Y_{\overline{k}}$-isomorphism $Y'_{\overline{k}} \to Z_{\overline{x}}$.

    Let $Z'' := Z \times_Y Y' = Z\times_{X \times Y} (X \times Y')$.
    Note that $Z''$ is normal, as the projection morphism to $Z$ is finite étale and $Z$ is normal.
    Moreover, since $X$ is normal geometrically connected and $Y'$ is normal connected, $X \times Y'$ is normal connected (and thus integral). Therefore,
    every connected component of $Z''$ is finite surjective over $X \times Y'$ and finite étale over $Z$.
    Let $Z_1'',\ldots,Z_t''$ be the connected components of $Z''$ and $I=\{1,\ldots,t\}$.
    For every $x \in \Sigma$, we have
    $Z''_{\overline{x}} = \coprod_{j\in I} Z''_{j,{\overline{x}}}$, and by choice of $\Sigma$,
    there exists an isomorphism $\{\overline{x}\} \times Y' \to Z_{\overline{x}}$ over $\{\overline{x}\} \times Y$, which induces a section of $Z''_{\overline{x}} \to \{\overline{x}\} \times Y'$ by the universal property of the fibre product. Since $Y'$ is geometrically irreducible, there exists $j(x) \in I$ such that this section restricts to a section of $Z''_{j(x),\overline{x}} \to \{\overline{x}\} \times Y'$.
    By the density of $\Sigma$ in $U(k)$, there exists $j_0 \in I$ such that the fibre of $j\colon \Sigma \to I$ over $j_0$ is dense in $U(k)$.
    Define $\Sigma' = j^{-1}(j_0)$ and $Z' := Z''_{j_0}$. Then $\Sigma'$ is a subset of $\Sigma$ such that $\Sigma'$ is dense in $U(k)$ and, for every $x\in \Sigma'$, the morphism $Z'_{\overline{x}} \to \{\overline{x}\} \times Y'$ is finite étale (since $Z \to X \times Y$ is étale over $U \times Y$) and has a section.

    Let $\overline{Y'} \to \overline{Y}$ be the normalization of $\overline{Y}$ in $Y'$.
    Let $\overline{Z'} \to X \times \overline{Y'}$ be the normalization of $X \times \overline{Y'}$ in $Z'$
    and note that this is étale over $U \times \overline{Y'}$.
    Finally, let $\overline{Z'} \to X' \to X$ be the Stein factorization of the composition $\overline{Z'} \to X \times \overline{Y'} \to X$
    and note that $X'$ is normal connected (thus integral).
    We have the following commutative diagram.

    \begin{equation}\label{Eq:VerticalRam:Diagram}
    \begin{gathered}
        \xymatrix{
            & \overline{Z'} \ar[dl] \ar[d] & Z' \ar[d] \ar[l] \ar[rr]^{\text{finite étale}} && Z \ar[d] \\
            X' \ar[dr] & X \times \overline{Y'} \ar[d] & X \times Y' \ar[d] \ar[l] \ar[rr]^{\text{finite étale}} && X \times Y \ar[d] \\
            & X \ar@{=}[r] & X \ar@{=}[rr] && X
        }
    \end{gathered}
    \end{equation}

    Let $Z' \to X' \times Y' = X' \times_{X} (X \times Y')$ be the natural morphism.
    We claim that $Z' \to X' \times Y'$ is an isomorphism and start by showing finiteness and surjectivity.
    By the definition of Stein factorization, the morphism $X' \to X$ is finite, and it is surjective since $Z' \to X$ is surjective and factorizes through it. 
    Since both $X' \times Y'$ and $Z'$ are finite over $X \times Y'$, the morphism $Z' \to X' \times Y'$ is also finite.

    Moreover, since both $Z' \to X \times Y'$ and $X' \times Y' \to X \times Y'$ are finite surjective, we have
    \[
        \dim (Z') = \dim (X\times Y') = \dim (X' \times Y')
    ,\]
    and the image of $Z'$ in $X' \times Y'$ is of dimension $\dim (X' \times Y')$.
    Since $X'$ is normal connected and $Y'$ is normal geometrically connected, $X' \times Y'$ is normal and connected,
    so
    $Z' \to X' \times Y'$ is dominant and thus (by its finiteness) surjective, as claimed.
   
    Next we prove that $Z' \to X' \times Y'$ is birational.
    Let $\overline{f}$ denote the morphism $\overline{Z'} \to X'$
    and $f = \overline{f}|_{Z'} \colon Z' \to X'$.
    Let $V \subseteq X'$ be a dense open such that the fibres of $\overline{f}^{-1}(V) \to V$ are normal.
    Since they are geometrically connected by the defining properties of Stein factorization, they are geometrically integral. By generic flatness and since a flat morphism of finite presentation is open, we may shrink $V$ further if necessary to assume that $f^{-1}(V) \to V$ is surjective.
    Thus, the fibres of $f^{-1}(V) \to V$ are nonempty and open in the fibres in $\overline{f}^{-1}(V) \to V$. Since the latter are normal and geometrically integral, the fibres of $f^{-1}(V) \to V$ are also normal and geometrically integral.
    The complement of $V$ in $X'$ maps to a proper closed subscheme of $X$. Let $U'$ denote the intersection of the complement of this closed subscheme with $U$.
    Then $U'$ is a dense open of $X$, contained in $U$, and since $X'_{U'} \subseteq V$, the fibres of $Z'_{U'} \to X'_{U'}$ are normal and geometrically connected.
    By generic étaleness and shrinking $U'$ further if necessary, we may and do assume that $X'_{U'} \to U'$ is étale.
    
    Let $x \in U'(k) \cap \Sigma'$ be a point. Then
    $X'_{\overline{x}}$ is finite étale over $\{\overline{x}\} = \Spec \overline{k}$,
    and thus it is the disjoint union $\coprod_{\overline{x}' \in X'_{\overline{x}}} \{\overline{x}'\}$,
    where each $\{\overline{x}'\}$ is a copy of $\Spec \overline{k}$.
    Let $g$ and $h$ denote the morphisms $Z' \to X' \times Y'$ and $Z' \to X \times Y'$, respectively, and let $g_{\overline{x}} \colon Z'_{\overline{x}} \to X'_{\overline{x}} \times Y'$
    and $h_{\overline{x}} \colon Z'_{\overline{x}} \to \{\overline{x}\} \times Y'$ denote the respective base changes.
    We then have the following commutative diagram.
    \[
        \xymatrix{
            \coprod_{\overline{x}' \in X'_{\overline{x}}} g_{\overline{x}}^{-1} (\{\overline{x}'\} \times Y') \ar@{=}[r] & Z'_{\overline{x}} \ar[d]_{g_{\overline{x}}} \ar[dr]^{h_{\overline{x}}}\\
            \coprod_{\overline{x}' \in X'_{\overline{x}}} (\{\overline{x}'\} \times Y') \ar@{=}[r] & X'_{\overline{x}} \times Y' \ar[r] \ar[d] & \{\overline{x} \} \times Y' \ar[d] \\
            \coprod_{\overline{x}' \in X'_{\overline{x}}} \{\overline{x}'\} \ar@{=}[r] & X'_{\overline{x}} \ar[r] & \{\overline{x} \}
        }
    \]
    Since $Z' \to X'$ has geometrically connected fibres over $U'$, the $g_{\overline{x}}^{-1} (\{\overline{x}'\} \times Y')$ are connected, so they are the connected components of $Z'_{\overline{x}}$.
    Moreover, by our choice of $\Sigma'$, the morphism $h_{\overline{x}}$ is finite étale and has a section.
    Thus, there exists a connected component $W$ of $Z'_{\overline{x}}$ such that $h_{\overline{x}}$ restricts to an isomorphism
    $h_{\overline{x}}|_W \colon W \to \{\overline{x}\} \times Y'$.
    Let $\overline{x}' \in X'_{\overline{x}}$ such that $g_{\overline{x}}^{-1} (\{\overline{x}'\} \times Y') = W$, i.e., $g_{\overline{x}}^{-1} (\{\overline{x}'\} \times Y') \to \{\overline{x}\} \times Y'$ is an isomorphism.
    Now choose a point $\overline{p} \in \{\overline{x}'\} \times Y'(\overline{k})$.
    Since both $\{\overline{x}'\} \times Y' \to \{\overline{x}\} \times Y'$ and
    $h_{\overline{x}}|_W$ are isomorphisms, there exists exactly one point
    in $W(\overline{k})$ over $\overline{p}$, i.e.,
    $|g_{\overline{x}}^{-1} (\overline{p})(\overline{k})| = 1$.
    Thus, $\overline{p}$ is a point of $(X' \times Y')(\overline{k})$ over $U'$ such that
    $g^{-1} (\overline{p})(\overline{k})$ has cardinality $1$.
    Now, since $g_{U'}\colon Z'_{U'} \to X'_{U'} \times Y'$ is finite étale, it follows that $g_{U'}$ is an isomorphism and we have shown birationality of $g$.
    As a finite surjective birational morphism of normal connected varieties, $g$ is an isomorphism, as claimed.

    We have now constructed the desired diagram.
    All of the claimed properties except for (\ref{Item:VerticalRam:Stein_is_ramified}) are true by construction or proved in the discussion above, so it only remains to verify that $X' \to X$ is a ramified cover.
    It is finite (as the finite part of the Stein factorization of $\overline{Z'} \to X$) and surjective (by the surjectivity of $Z'\to X$).
    Thus, it remains to show that $X' \to X$ is ramified. For a contradiction, assume that it is unramified (and therefore étale).
    Then $Z' \cong X' \times Y' \to X \times Y' \to X \times Y$ and thus $Z' \to Z \to X \times Y$ would be étale. Since $Z' \to Z$ is étale surjective, $\pi$ would be étale by the cancellation property of étale morphisms, contradiction.
\end{proof}

\section{Arithmetic refinements}

If $S$ is
a normal noetherian integral scheme with function field $k$ of characteristic zero
and $\mathcal{X}$ is an arithmetic scheme over~$S$, we ``test'' the weak Hilbert property of $\mathcal{X}$ over $S$ by considering finite families of ramified covers
$(\pi_i \colon Z_i \to \mathcal{X}_k)_{i=1}^n$. When doing this, it might be beneficial to replace a given cover $\pi_i$ by taking finitely many covers $(f_{ij} \colon W_{ij} \to Z_i)_{j\in J_i}$ such that
\[
    \mathcal{X}(S)^{(1)} \setminus \bigcup_{i=1}^n \pi_i(Z_i(k))
    = \mathcal{X}(S)^{(1)} \setminus \bigcup_{i=1}^n\bigcup_{j \in J_i} (\pi_i \circ f_{ij})(W_{ij}(k))
.\]
This is achieved, for example, if $Z_i(k) = \bigcup_{j \in J_i} f_{ij}(W_{ij}(k))$, in which case the collection $(f_{ij} \colon W_{ij} \to Z_i)_{j\in J_i}$ is called an \emph{arithmetic refinement} of the variety $Z_i$ in \cite[§3.2]{CDJLZ}.
Arithmetic refinements of smooth proper varieties are constructed in \emph{loc.\ cit.}\ as an application of the Chevalley--Weil Theorem, which we can not apply in order to lift $k$-rational points, since the varieties in our situation are not proper.
Nevertheless (using the same proof of the Chevalley--Weil Theorem), it is possible for us to lift near-integral points (cf. Theorem~\ref{Thm:ChevalleyWeilIntegral}), which leads us to the following definition.

\begin{definition}
    Let $S$ be
    a normal noetherian integral scheme with function field $k$ of characteristic zero and let $\mathcal{Z}$ be an arithmetic scheme over $S$. An \emph{arithmetic refinement} of $\mathcal{Z}$ is a finite family
    $(f_j \colon W_j \to \mathcal{Z}_k )_{j \in J}$ of covers with $W_j$ normal and integral over $k$ such that
    $\mathcal{Z}(S)^{(1)} \subseteq \bigcup_{j\in J} \psi_j(W_j(k))$.
\end{definition}

In order to construct suitable arithmetic refinements, we will use the following Chevalley--Weil type lifting theorem for near-integral points. The proof of this theorem is the same as the proof of the Chevalley--Weil Theorem for finitely generated fields of characteristic $0$ given in \cite[Theorem~3.8]{CDJLZ}.

\begin{theorem}\label{Thm:ChevalleyWeilIntegral}
    Let $S$ be
    a regular integral arithmetic scheme
    and $f \colon \mathcal{X} \to \mathcal{Y}$ a finite étale morphism of integral finite type schemes over $S$.
    Then there exists a finite étale morphism $S' \to S$
    of regular integral arithmetic schemes
    such that
    $\mathcal{Y}(S)^{(1)} \subseteq f(\mathcal{X}(S')^{(1)})$.
\end{theorem}
\begin{proof}
    Let $y \in \mathcal{Y}(S)^{(1)}$ be a near-integral point, so there exists a dense open subscheme $U_y \subseteq S$ with complement of codimension at least two such that $y$ is induced by a morphism $U_y \to \mathcal{Y}$.
    Let $V_y$ be a connected component of $U_y \times_{\mathcal{Y}} \mathcal{X} \to U_y$
    and note that $V_y \to U_y$
    is finite étale surjective and of degree at most $\deg(f)$.
    Let $\overline{V_y} \to S$ be the normalization of $S$ in the function field of $V_y$.
    
    We claim that $\overline{V_y} \setminus V_y$ is of codimension at least $2$ in $\overline{V_y}$.
    Indeed, let $D \subseteq \overline{V_y} \setminus V_y$ be an irreducible component. Note that $D$ maps finitely to $S\setminus U_y$, so we have $\dim D = \dim (S \setminus U_y) \leq \dim S - 2$.
    Since $\dim S = \dim \overline{V_y}$ by the finiteness of $\overline{V_y} \to S$, this proves the claim.
    By purity of the branch locus \cite[Théorème~X.3.1]{SGA1}, it follows that $\overline{V_y} \to S$ is étale.

    Since there are only finitely many isomorphism classes of finite étale covers of $S$ with bounded degree \cite[Theorem~2.9]{HaradaHiranouchi}, the set of isomorphism classes of the morphisms $\overline{V_y} \to S$ with $y \in \mathcal{Y}(S)^{(1)}$ is finite.
    Thus, there exists a finite étale morphism $S' \to S$ of
    integral noetherian
    schemes
    such that, for every $y \in \mathcal{Y}(S)^{(1)}$,
    the morphism $S' \to S$ factors over $\overline{V_y}$,
    i.e., we have $\mathcal{Y}(S)^{(1)} \subseteq f(\mathcal{X}(S')^{(1)})$.
    Note that the composition $S' \to S \to \Spec \mathbb{Z}$ is dominant of finite type
    and that $S'$ is regular since $S' \to S$ is finite étale and $S$ is regular.
    Thus, $S'$ is a regular integral arithmetic scheme, as required.
\end{proof}

\begin{corollary}\label{Cor:CWT_NearIntegralToRational}
    Let $S$ be a regular integral arithmetic scheme, let $\mathcal{Y}$ be an integral finite type scheme over $S$, and let $f \colon X \to \mathcal{Y}_{K(S)}$ be a finite étale morphism of integral $K(S)$-schemes.
    Then there exists a finite field extension $L/k$ such that $\mathcal{Y}(S)^{(1)} \subseteq f(X(L))$.
\end{corollary}
\begin{proof}
    By standard spreading out arguments, there exists a dense open $S' \subseteq S$,
    an integral finite type $S'$-scheme $\mathcal{X}$ with $\mathcal{X}_{K(S')} = X$,
    and a finite étale morphism $\mathcal{X} \to \mathcal{Y}_{S'}$ extending $f$.
    By Theorem~\ref{Thm:ChevalleyWeilIntegral}, there exists a finite étale (in particular, finite surjective) morphism $S'' \to S'$ of integral arithmetic schemes such that
    $\mathcal{Y}_{S'}(S')^{(1)} \subseteq f(\mathcal{X}(S'')^{(1)})$.
    Since $\mathcal{Y}(S)^{(1)} \subseteq \mathcal{Y}_{S'}(S')^{(1)}$ (as $S' \subseteq S$ is a dense open), the claim follows by taking $L=K(S'')$.
\end{proof}

We can now prove a generalization of \cite[Theorem~3.9]{CDJLZ} that allows us to replace a vertically ramified cover by a suitable arithmetic refinement (of an arithmetic model) when testing for the weak Hilbert property of a product of arithmetic schemes.

\begin{theorem}\label{Thm:ArithmeticRefinement}
    Let $S$ be
   a regular integral arithmetic scheme
    with function field $k$, let $\mathcal{X}$ and~$\mathcal{Y}$ be quasi-projective arithmetic schemes over $S$,
    and define $X := \mathcal{X}_k$ and $Y := \mathcal{Y}_k$.
    Let $\pi \colon Z \to X \times Y$ be a cover which is vertically ramified over $X$
    and assume that $Z(k)$ is dense in $Z$, so that in particular $X(k)$ and $Y(k)$ are dense in $X$ and $Y$, respectively.
    Let $Y \to \overline{Y}$ be an open immersion into a normal projective integral variety over $k$ and let $\overline{Z} \to X \times \overline{Y}$ be the normalization of $X \times \overline{Y}$ in $Z$.
    Then there exists a finite family of covers $(f_j \colon W_j \to Z)_{j\in J}$ such that
    \[
        (\mathcal{X} \times_S \mathcal{Y})(S)^{(1)}
        \setminus \pi(Z(k))
        = (\mathcal{X} \times_S \mathcal{Y})(S)^{(1)}
        \setminus \bigcup_{j\in J} (\pi \circ f_j)(W_j(k))
    \]
    and, for every $j\in J$, the Stein factorization of $\overline{W_j} \to X\times \overline{Y} \to X$ is ramified over $X$, where $\overline{W_j} \to \overline{Z}$ is the normalization of $\overline{Z}$ in $W_j$.
\end{theorem}
\begin{proof}
    Let $\overline{Z} \to T \to X$ be the Stein factorization of the composition $\overline{Z} \to X \times \overline{Y} \to X$ and note that $T$ is a normal geometrically integral $k$-variety. 
    If $T \to X$ is ramified,
    the family consisting of the identity map $f = \id \colon W = Z \to Z$ proves the claim,
    so we may assume that $T \to X$ is étale.

    The morphism $\overline{Z} \to X \times \overline{Y}$ factors over a morphism $\overline{Z} \to T \times \overline{Y}$.
    Since $\overline{Z}$ is the normalization of $X \times \overline{Y}$ in $Z$, we have that $\overline{Z}$ is also the normalization of $T\times \overline{Y}$ in $Z$.
    Moreover, $\overline{Z} \to T \times \overline{Y}$ is proper and quasi-finite by the finiteness of both $\overline{Z} \to X \times \overline{Y}$ and $T \times \overline{Y} \to X \times \overline{Y}$, so it is finite surjective.
    Let $\tilde{Z} = \overline{Z} \times_{T \times \overline{Y}} (T \times Y)$.
    Note that $\tilde{Z} = \overline{Z} \times_{X \times \overline{Y}} (X \times Y)$, so there is a natural morphism $Z \to \tilde{Z}$, which is proper (and in particular closed) by the properness of $Z \to X \times Y$ \cite[Tag~01W6]{stacks-project}.
    Moreover, since $Y \to \overline{Y}$ (and thus $\tilde{Z} \to \overline{Z}$) and $Z \to \overline{Z}$ are open immersions, it follows that $Z \to \tilde{Z}$ is an open immersion.
    By the connectedness of $\tilde{Z}$, we conclude that $Z \to \tilde{Z}$ is a surjective open immersion and therefore an isomorphism.
    In particular, $Z \to T\times Y$ is finite surjective.
    As $T\times Y \to X \times Y$ is étale and $Z \to X \times Y$ is ramified, $Z \to T \times Y$ is ramified. Moreover, by the étaleness of $Z$ over $U \times Y$ and the cancellation law for étale morphisms, $Z$ is étale over $T_U \times Y$, i.e., $Z \to T \times Y$ is a cover vertically ramified over $T$.

    Since the composition $\overline{Z} \to T \times \overline{Y} \to T$ has geometrically connected fibres by the defining properties of Stein factorization and $Z(k)$ is dense, 
    Lemma~\ref{Lem:VerticalRamification} gives us a diagram
    \begin{equation}\label{Eq:DiagramRefinements}
    \begin{gathered}
    \xymatrix{
        & \overline{Z'} \ar[dl] \ar[d] & Z' \cong T' \times Y' \ar[l] \ar[d] \ar[r] & Z \ar[d] \\
        T' \ar[dr] & T \times \overline{Y'} \ar[d] & T \times Y' \ar[l] \ar[d] \ar[r] & T \times Y \ar[d] \\
        & T \ar@{=}[r] & T \ar@{=}[r] & T
     }
    \end{gathered}
    \end{equation}
    of varieties over $k$ such that $T' \to T$ is a ramified cover, $Y' \to Y$ and $Z'\to Z$ are finite étale, $\overline{Y'}$ is the normalization of $\overline{Y}$ in $Y'$, and $\overline{Z'}$ is the normalization of $T \times \overline{Y'}$ in $Z'$.
    Let $\psi$ denote the morphism $ Z' \to Z$.

    By standard spreading out arguments,
    there exist a regular dense open subscheme $S' \subseteq S$
    and an integral finite type separated $S'$-scheme $\mathcal{Z}$ whose generic fibre over $S'$ is $Z$, and
    such that $\psi$ and $\pi$ extend to a finite étale (hence surjective, since $\mathcal{Z}$ is connected) morphism
    $\mathcal{Z'} \to \mathcal{Z}$
    and a finite surjective morphism
    $\mathcal{Z} \to (\mathcal{X} \times \mathcal{Y})_{S'}$,
    respectively.
    If $x\in (\mathcal{X} \times_S \mathcal{Y})(S)^{(1)}$ is a near-$S$-integral point that lifts along $\pi$ to a $k$-rational point of $Z$, then, by the valuative criterion of properness, we have $x \in \pi(\mathcal{Z}(S')^{(1)})$.

    By Corollary~\ref{Cor:CWT_NearIntegralToRational}, there is a finite field extension $L/k$ such that
    $\mathcal{Z}(S')^{(1)} \subseteq \psi(Z'(L))$.
    In the next steps, we follow the arguments of \cite[Theorem~3.9]{CDJLZ} to construct an arithmetic refinement of $\mathcal{Z}$.
    Consider the induced morphism of Weil restriction of scalars (cf. \cite{BoschNeronModels})
    \[
        R_{L/k}(\psi) \colon R_{L/k} (Z'_L) \to  R_{L/k} (Z_L)
    ,\]
    which is finite étale (as $\psi$ is).
    Let $\Delta \colon Z \to R_{L/k}(Z_L)$ be the diagonal morphism
    and define
    $W = Z \times_{\Delta, R_{L/k}(Z_L)} R_{L/k} (Z'_L)$,
    so we have the following cartesian diagram.
    \[
        \xymatrixcolsep{5pc}\xymatrix{
            W \ar[d] \ar[r] & Z \ar[d]^{\Delta} \\
             R_{L/k}(Z'_L) \ar[r]^{R_{L/k}(\psi)} &  R_{L/k}(Z_L)
        }
    \]
    Since $\psi$ (and thus $R_{L/k}(\psi)$) is finite étale, $W$ is finite étale over $Z$, so it is normal.
    Let $W_j, j \in J,$ be the connected components of $W$ that satisfy $W_j(k) \neq \emptyset$
    and note that these are normal, geometrically connected (thus geometrically integral),
    and finite étale over $Z$.
    
    The map of $L$-points $Z'(L) \to Z(L)$ induced by $\psi$ is identified by $R_{L/k}$ with
    $R_{L/k}(Z'_L)(k) \to  R_{L/k}(Z_L)(k)$
    and the map $Z(k) \to R_{L/k}(Z_L)(k) \cong Z(L)$ induced by $\Delta$ is the inclusion map.
    Thus, for every point $x \in Z(k)$ that lifts along $\psi$ to $Z'(L)$,
    there exists a point in $R_{L/k}(Z'_L)(k)$ lying over $\Delta(x)$,
    so $x$ lifts to $W(k)$.
    In particular, $\mathcal{Z}(S)^{(1)}$ is contained in the image of $W(k)$ in $Z(k)$, i.e.,
    the family $(f_j \colon W_j \to Z)_{j\in J}$ is an arithmetic refinement of $\mathcal{Z}$ and we have
    \begin{align*}
        (\mathcal{X} \times_S \mathcal{Y})(S)^{(1)}
        \setminus \pi(Z(k))
        & = (\mathcal{X} \times_S \mathcal{Y})(S)^{(1)}
        \setminus \pi(\mathcal{Z}(S')^{(1)}) \\
        & = (\mathcal{X} \times_S \mathcal{Y})(S)^{(1)}
        \setminus \bigcup_{j\in J} (\pi \circ f_j)(W_j(k))
    .\end{align*}
    To finish the proof, it remains to show that the Stein factorization of the composed morphism $\overline{W_j} \to X \times \overline{Y} \to X$ is ramified.
    Consider the composed morphism $W \to Z \to T$ and, for every $j\in J$,
    let $T_j$ denote the normalization of $T$ in $W_j$.
    Note that $T_j$ is normal and integral by \cite[Tag~035L]{stacks-project} (in particular, $T_j \to T$ is a cover)
    and, since $T \to X$ is finite, $T_j$ is also the normalization of $X$ in $W_j$.
    
    The base change of $W = Z \times_{\Delta, R_{L/k}(Z_L)} R_{L/k}(Z'_L)$ to the algebraic closure $\overline{k}$ is given by
    $W_{\overline{k}} \cong Z'_{\overline{k}} \times_{Z_{\overline{k}}} \ldots \times_{Z_{\overline{k}}} Z'_{\overline{k}}$,
    so we have a finite étale morphism $W_{\overline{k}} \to Z'_{\overline{k}}$ over $Z_{\overline{k}}$, and similarly a finite étale morphism $W_{j,\overline{k}} \to Z'_{\overline{k}}$ over $Z_{\overline{k}}$ for every $j \in J$.
    
    Since normalization commutes with flat base change (in particular, base change from $k$ to $\overline{k}$),
    $T_{j,\overline{k}}$ is the normalization of $T_{\overline{k}}$ in $W_{j,\overline{k}}$.
    By the universal mapping property of normalization \cite[Tag~035I]{stacks-project},
    there is a finite surjective morphism $T_{j,\overline{k}} \to T'_{\overline{k}}$ over $T_{\overline{k}}$.
    Étaleness can be checked after base change to $\overline{k}$, so $T'_{\overline{k}} \to T_{\overline{k}}$ is ramified (as $T' \to T$ is).
    If $T_{j,\overline{k}} \to T'_{\overline{k}}$ is ramified, then $T_{j,\overline{k}} \to T_{\overline{k}}$ is ramified by \cite[Tag~02GG]{stacks-project},
    and if $T_{j,\overline{k}} \to T'_{\overline{k}}$ is étale, then $T_{j,\overline{k}} \to T_{\overline{k}}$ is ramified by \cite[Tag~02K6]{stacks-project}.
    We conclude that $T_{j,\overline{k}} \to T_{\overline{k}}$ is ramified, which implies that
    This implies that $T_j \to T$ and thus the composition $T_j \to T \to X$ is ramified.

    Let $\overline{W_j} \to \overline{T_j} \to X$ be the Stein factorization of $\overline{W_j} \to X$, i.e., $\overline{T_j}$ is the normalization of $X$ in $\overline{W_j}$.
    Since $W_j \subseteq \overline{W_j}$ is an open subscheme, we have an open immersion
    $T_j \hookrightarrow \overline{T_j}$ over $X$.
    Thus, as $T_j \to X$ is ramified, we conclude that $\overline{T_j} \to X$ is ramified, as required.
\end{proof}

\section{The weak Hilbert property}

Recall that an arithmetic scheme $\mathcal{X}$ over a normal noetherian integral scheme $S$  with $K(S)$ of characteristic $0$ has the weak Hilbert property over $S$
if $\mathcal{X}(S)^{(1)}$ is not strongly thin (Definition~\ref{Def:WHP}).
The usual Hilbert property is sometimes also formulated using the integrality of many fibres.
For the sake of completeness, we include the following related observation for the weak Hilbert property.

\begin{proposition}\label{Prop:WHP_IntegralFibres}
    Let $X$ be a normal integral variety over a field $k$ of characteristic $0$
    and $\Sigma \subseteq X(k)$ a subset.
    Let $Z \to X$ be a Galois ramified covering with no nontrivial étale subcovers.
    If $\Sigma$ is not strongly thin in $X$, there exists a dense subset of points $x \in \Sigma$ such that $Z_x$ is integral.
\end{proposition}
\begin{proof}
    Let $G$ be the Galois group of the cover $Z \to X$.
    By assumption, for every normal subgroup $H \subsetneq G$, the cover given by the quotient $Z/H \to X$ is ramified.
    Since $\Sigma \subseteq X(k)$ is not strongly thin and by generic étaleness, we find a dense set of points $x\in \Sigma$ such that, for all normal subgroups $H\subsetneq G$, the cover $Z/H \to X$ is étale over $x$ and $(Z/H)_x(k) = \emptyset$.
    We show that, for all of those $x$, the fibre $Z_x$ is integral.
    Since $G$ acts transitively on $Z_x$, we have $Z_x = \Spec \prod_{i=1}^n L$ with $L$ a Galois extension of $k$ and $n[L:k] = \deg (Z \to X) = \# G$.
    Let $H$ be the Galois group of $L/k$. If $H=G$, we have $n=1$ and thus irreducibility of $Z_x$. If $H \neq G$, then $(Z/H)_x$ has a $k$-point, contradicting the choice of $x$.
\end{proof}

\begin{remark}
    Note that, if $Z \to X$ is a ramified cover with no nontrivial étale subcovers and $Z' \to Z \to X$ its Galois closure, it is not necessarily true that $Z' \to X$ has no nontrivial étale subcovers.
    For example, let $A = B = \mathbb{P}^1_{\mathbb{Q}}$ and take a degree $3$ cover $B \overset{\pi}{\to} A$, ramified over three $\mathbb{Q}$-points of $A$ and not Galois.
    Let $C \to B \to A$ be the Galois closure of $\pi$, so that $C \to A$ is of degree $6$, and let $C \to D \to A$ be a subcover with $D \to A$ of degree $2$.
    Note that $D \to A$ can only be ramified over branch points of $\pi$ and thus has at most three branch points. It now follows from the Riemann--Hurwitz formula that $D$ is of genus $0$, so $D_{\overline{\mathbb{Q}}} \to A_{\overline{\mathbb{Q}}}$ has precisely two branch points.
    In particular, the étale locus $X \subseteq A$ of $D \to A$ contains at least one branch point of $\pi$. Define $Z := B_X$ and $Z' := C_X$.
    Then $Z \to X$ is a ramified cover (with no nontrivial subcovers as its degree is prime), $Z' \to Z \to X$ is its Galois closure, and $Z' \to D_X \to X$ is a nontrivial étale subcover.
\end{remark}

In the remainder of this section, we are interested in subsets of near-integral points that are not strongly thin, and whether this property persists along certain morphisms.
To do so, we use the following definition.

\begin{definition}\label{Def:WeaklyHilbertianPair}
	Let $X$ be a normal integral variety over a field $k$ of characteristic $0$ and $\Sigma \subseteq X(k)$ a subset.
	We say that $(X, \Sigma)$ is a \emph{weakly Hilbertian pair (over $k$)}, or \emph{weakly Hilbertian}, if $\Sigma$ is not strongly thin in $X$.
\end{definition}
In particular,
if $S$ is a normal noetherian integral scheme with function field $k$ of characteristic~$0$,
then an arithmetic scheme $\mathcal{X}$ over $S$
has the weak Hilbert property over $S$
if and only if the pair $(\mathcal{X}_k, \mathcal{X}(S)^{(1)})$ is weakly Hilbertian.

In \cite[Definition~1.2]{JavanpeykarNefTangentBundle}, a pair $(X, \Omega)$ (with $X$ a normal projective variety over a field $k$ of characteristic $0$ and $\Omega\subseteq X(k)$) is called Hilbertian if $\Omega$ is not strongly thin in $X_L$ for any finite extension $L/k$.
It can be seen from \cite[Lemma~2.1]{BSFP23BaseChange} that this is in fact equivalent to $\Omega$ not being strongly thin in $X$. Thus, Definition~\ref{Def:WeaklyHilbertianPair} is a generalization of \cite[Definition~1.2]{JavanpeykarNefTangentBundle}.

In what follows, we combine the results of the previous sections to deduce our main theorem (Theorem~\ref{Thm:ProductsWHP}),
following the proof of \cite[Theorem~1.9]{CDJLZ}.
After the proof of our product theorem, we show some persistence properties of the weak Hilbert property of arithmetic schemes, which we deduce from slight generalizations of results obtained in \cite[§3.1 and §3.6]{CDJLZ}. For the most part, the proofs of these generalizations are already contained in the proofs of \emph{loc.\ cit.}, up to some minor changes.
Therefore, we repeat the parts of the proofs given in \cite{CDJLZ} where we make changes and those necessary to introduce notation, and omit the parts that are the same for our situation.

\subsection{Product theorem}

To prove the announced Theorem~\ref{Thm:ProductsWHP}, we show the following more general statement, which is close in spirit to \cite[Theorem~1.1]{BSFP14}:

\begin{theorem}\label{Thm:FibrationProducts}
    Let $S$ be a regular integral arithmetic scheme
    with function field $k$,
    let $\mathcal{X}, \mathcal{Y}$ be quasi-projective arithmetic schemes over $S$,
    and define $X=\mathcal{X}_k, Y = \mathcal{Y}_k$.
    Let $p \colon X \times Y \to X$ denote the projection morphism.
    Let $\Sigma \subseteq (\mathcal{X} \times_S \mathcal{Y})(S)^{(1)}$ and define
	$\Sigma_X = p(\Sigma) \subseteq \mathcal{X}(S)^{(1)}$.
    For $x \in X(k)$, define $\Sigma_x = \Sigma\cap (\{x\}\times Y)$.
    If $(X, \Sigma_X)$ is weakly Hilbertian and $(\{x\}\times Y, \Sigma_x)$ is weakly Hilbertian for every $x\in \Sigma_X$,
	then $(X \times Y, \Sigma)$ is weakly Hilbertian.
\end{theorem}
\begin{proof}
    Note first that $X$ and $Y$ have a dense set of $k$-rational points by the assumptions, so they are geometrically irreducible by \cite[Tag~0G69]{stacks-project}, and hence $X \times Y$ is a normal (geometrically) integral $k$-variety.
    Let $\left(\pi_i \colon Z_i \to X\times Y\right)_{i=1}^n$ be a finite family of ramified covers with $Z_i$ normal geometrically integral varieties over $k$.
    Without loss of generality, we may assume that $Z_i(k)$ is dense in $Z_i$ for every $i$ (as otherwise the image $\pi_i(Z_i(k))$ is not dense in $X \times Y$).
    Let $Y \to \overline{Y}$ be an open immersion with $\overline{Y}$ a normal projective geometrically integral variety over $k$,
    and let $\overline{Z_i} \to X \times \overline{Y}$ be the normalization of $X \times \overline{Y}$ in $Z_i$.
    Let $\psi_i \colon T_i \to X$ denote the Stein factorization of $\overline{Z_i} \to X$.
    
    Let $i \in \{1,\ldots, n\}$ be an index such that
    $Z_i \to X\times Y$ is vertically ramified over $X$.
    By Theorem~\ref{Thm:ArithmeticRefinement}, there exists a finite family of covers
    $(f_{ij} \colon W_{ij} \to Z_i )_{j\in J_i}$ such that
    \[
        \Sigma
        \setminus \pi_i(Z_i(k))
        = \Sigma
        \setminus \bigcup_{j\in J_i} (\pi_i \circ f_{ij})(W_{ij}(k))
    \]
    and, for every $j\in J_i$, the Stein factorization of $\overline{W_{ij}} \to X\times \overline{Y} \to X$ is ramified over $X$, where
    $\overline{W_{ij}} \to \overline{Z_i}$ is the normalization of $\overline{Z_i}$ in $W_{ij}$.
    Thus, we may replace $Z_i \to X\times Y$ by the finitely many covers
    $W_{ij} \to X\times Y$ to assume that $\psi_i$ is ramified.
    Up to reindexing, we may now assume that there exists an integer $m\in \{1,\ldots,n\}$ such that
    \begin{enumerate}
        \item $\psi_i$ is ramified for all $i=1,\ldots,m$,
        \item $\psi_i$ is étale and the branch locus of $\pi_i$ dominates $X$ for all $i=m+1,\ldots,n$.
    \end{enumerate}

    Let
    \[
        \Omega := \Sigma_X \setminus \bigcup_{i=1}^m \psi_i(T_i(k))
    ,\]
    which is dense in $X$ since $\Sigma_X$ is not strongly thin.
    To conclude the proof, we show that
    \[
        \Psi
        := \bigcup_{x\in \Omega} \left( \Sigma_x \setminus \bigcup_{i=m+1}^n\pi_{i,x}(Z_{i,x}(k)) \right)
        \subseteq \Sigma \setminus \bigcup_{i=1}^n\pi_{i}(Z_{i}(k))
    \]
    is dense in $X \times Y$.

    Let $V \subseteq X$ be a dense open subscheme such that, for every $x$ in $V(k)$ and every $i\in \{m+1,\ldots,n\}$, the fibre
    $(\overline{Z_i})_x$
    is normal.
    Recall that, for $i\in \{m+1,\ldots,n\}$, the ramification locus $D_i \subseteq X \times Y$ of $\pi_i$ dominates $X$. Thus, the image of $D_i$ in $X$ is a dense constructible set and therefore contains a dense open subset $U_i \subseteq X$. Let $U = V \cap \bigcap_{i=m+1}^n U_i$.
    
    Let $\pi_{i,T_i} \colon Z_{i,T_i} \to T_i \times Y$ denote the pullback of $\pi_i$ along $\psi_{i,Y} \colon T_i \times Y \to X \times Y$. Since $\psi_{i,Y}$ is étale, the ramification locus $D_{i,T_i} \subseteq T_i \times Y$ of $\pi_{i,T_i}$ is given by the pullback of $D_i$ along $\psi_{i,Y}$.

    By the density of $\Omega \cap U$ in $X$, we obtain the desired density of $\Psi$ by showing that, for every $x \in \Omega \cap U$, the set
    $\Sigma_x \setminus \bigcup_{i=m+1}^n \pi_{i,x}(Z_{i,x}(k))$
    is dense in $\{x\} \times Y \cong Y$.
    Since $\Sigma_x$ is not strongly thin, this follows after showing that, for $m+ 1 \leq i \leq n$ and $x \in U(k)$, and for every irreducible component $Z'$ of $Z_{i,x}$, the restriction of $\pi_{i,x}$ to $Z'$ is a ramified cover.

    Let $i\in \{m+1, \ldots, n\}$ and fix $x\in U(k)$.
    Let $\psi_i^{-1}(x) = t_1 \sqcup \ldots \sqcup t_r$ be the decomposition into isolated points (recall that $\psi_i$ is étale).
    By its normality, $(\overline{Z_i})_x = (\overline{Z_i})_{t_1} \sqcup \ldots \sqcup (\overline{Z_i})_{t_r}$ decomposes as a disjoint union of normal varieties. Since the fibres of $\overline{Z_i} \to T_i$ are geometrically connected, the varieties $(\overline{Z_i})_{t_j}$ are connected (and thus integral).
    Let $Z'$ be an irreducible component of $Z_{i,x}$. Then $Z' = Z_{i,t_j} \subseteq (\overline{Z_i})_{t_j}$ is an open subscheme for some point $t_j \in T_{i,x}$.
    Since $D_{i,T_i}$ is given by the pullback of $D_i$ along $\psi_{i,Y}$ and $x\in U(k)$ lies in the image of $D_i$ inside $X$, the universal property of fibre products gives us a point of $D_{i,T_i}$ lying over $t_j$.
    This implies that $Z' = Z_{i,t_j} \to \{t_j\} \times Y$ and thus $Z' \to \{x\}\times Y$ ramifies, as claimed.
\end{proof}

We obtain Theorem~\ref{Thm:ProductsWHP} as an immediate corollary.

\begin{proof}[Proof of Theorem~\ref{Thm:ProductsWHP}]
    Apply Theorem~\ref{Thm:FibrationProducts} to the set of near-integral points
    $\Sigma := (\mathcal{X} \times_S \mathcal{Y})(S)^{(1)} = (\mathcal{X}(S)^{(1)}) \times (\mathcal{Y}(S)^{(1)})$.
\end{proof}

\subsection{Birational invariance}

The weak Hilbert property is a birational invariant among smooth proper geometrically connected varieties over a finitely generated field of characteristic $0$, as is proved in \cite[Proposition~3.1]{CDJLZ}.
Adapting the arguments of the proof of \emph{loc.\ cit.}, we show that Hilbertianity of a pair is invariant under proper birational morphisms in the following sense.

\begin{lemma}\label{Lem:HPairsBirational}
	Let $\pi \colon X' \to X$ be a proper birational morphism of smooth integral varieties over a field $k$ of characteristic $0$, let $\Sigma_X \subseteq X(k)$ be a subset, and define $\Sigma_{X'}$ to be the preimage of $\Sigma_X$ under the induced morphism $X'(k) \to X(k)$.
    Then the pair $(X', \Sigma_{X'})$ is weakly Hilbertian if and only if $(X, \Sigma_{X})$ is weakly Hilbertian.
\end{lemma}
\begin{proof}
    Let $\Cov(X)$ (resp. $\Cov(X')$) denote the class of covers $Y \to X$ (resp. $Y \to X'$) with $Y$ a normal integral $k$-variety.
    Define the maps $N' \colon \Cov(X) \to \Cov(X')$ and $N \colon \Cov(X') \to \Cov(X)$ as follows:
    
    If $Y' \to X'$ is a cover, let $N(Y') \to X$ be the normalization of $X$ in $Y'$.
    If $Y \to X$ is a cover, note that the generic fibre of $Y\times_{X} X' \to X'$ is isomorphic to the generic fibre of $Y \to X$ (since $X' \to X$ is generically an isomorphism) and thus irreducible. Therefore, $Y\times_{X} X'$ has exactly one irreducible component $Y'_0$ dominating $X'$. The morphism $Y'_0 \to Y$ is birational and $Y'_0 \to X$ is finite surjective since $Y'_0 \to Y\times_{X} X'$ is a closed immersion and thus finite.
    We now let $N'(Y) \to X'$ be the normalization of $X'$ in $Y'_0$.
    Note that both $N$ and $N'$ preserve the degree of a cover by the birationality of $X'$ and $X$.

    We claim that $N$ and $N'$ are inverse to each other.
    To see this, 
    let $Y \to X$ be a cover with $Y$ normal and define $Y'_0$ as above.
    By the universal property of normalization, the inclusion $Y'_0 \to Y \times_X X'$ factors as $Y'_0 \to N'(Y) \to Y \times_X X'$ over $X'$. In particular, this induces a factorization $N'(Y) \to Y \to X$ of the composed morphism $N'(Y) \to X' \to X$.
    Therefore, again by the universal property of normalization,
    we have a morphism $N(N'(Y)) \to Y$ over $X$, which is finite (by the finiteness of $N(N'(Y)) \to X$) and birational (by the degrees over $X$ being equal), and thus an isomorphism by Zariski's main theorem.
    Conversely, if $Y' \to X'$ is a cover with $Y'$ normal,
    we obtain a finite birational morphism $N'(N(Y')) \to Y'$ over $X'$, which again is an isomorphism by Zariski's main theorem, as required.

    Next we show that $N$ and $N'$ send étale covers to étale covers.
    If $Y \to X$ is an étale cover, then $Y \times_{X} X'$ is étale over $X'$ and thus normal.
    Hence, the irreducible component $Y'_0$ dominating $X'$ is normal and étale over $X'$ and we have $N'(Y) = Y'_0$.
    
    For the converse, let $Y' \to X'$ be an étale cover and consider the cover
    $N(Y') \to X$.
    Since $X' \to X$ is proper birational, by \cite[Corollary~4.4.3]{LiuAGAC}, there exists a dense open $U \subseteq X$ such that
    $U' := \pi^{-1}(U) \to U$ is an isomorphism and the codimension of $X\setminus U$ in $X$ is at least~$2$. The pullback $N(Y')\times_X U \to U$ is étale, as its pullback along the isomorphism $U' \to U$ is $Y' \times_{X'} U' \to U'$ (which is étale).
    Since $X \setminus U$ is of codimension at least $2$, it follows from purity of the branch locus \cite[Théorème~X.3.1]{SGA1} that $N(Y') \to X$ is étale, as claimed.
    We conclude that $N$ and $N'$ are inverse to each other and send étale covers to étale covers and thus ramified covers to ramified covers.

    To finish the proof, since a subset is strongly thin if and only if its intersection with some dense open subscheme is strongly thin, we may assume that $\Sigma_X \subseteq U(k)$ (and thus $\Sigma_{X'} \subseteq U'(k)$).
    Given a finite family of ramified covers $(\pi_i \colon Y_i \to X)_{i=1}^n$
    (resp. $(\pi_i' \colon Y_i' \to X)_{i=1}^n$), we define $Y_i' = N'(Y_i)$ (resp. $Y_i = N(Y_i')$)
    to obtain a finite family of ramified covers $(\pi_i' \colon Y_i' \to X)_{i=1}^n$
    (resp. $(\pi_i \colon Y_i \to X)_{i=1}^n$) such that $\pi$ induces a bijection
    \[
        \Sigma_{X'} \setminus \bigcup_{i=1}^n \pi_i'(Y_i'(k)) \cong
        \Sigma_X \setminus \bigcup_{i=1}^n \pi_i(Y_i(k))
    .\]
    In particular, the left side is dense in $X'$ if and only if the right side is dense in $X$, which finishes the proof.
\end{proof}

In order to derive a birationality theorem for the weak Hilbert property for arithmetic schemes from Lemma~\ref{Lem:HPairsBirational}, we have to ensure that near-integral points lift along the considered birational morphism, which leads to the following statement.

\begin{proposition}\label{Prop:WHPBirational}
    Let $\mathcal{X}, \mathcal{X}'$ be arithmetic schemes over a normal integral arithmetic scheme $S$. Assume that $\mathcal{X}_{K(S)}$ and $\mathcal{X}'_{K(S)}$ are smooth and that there exists a proper birational morphism $\pi \colon \mathcal{X}' \to \mathcal{X}$ over $S$.
    Then $\mathcal{X}$ has the weak Hilbert property over $S$ if and only if $\mathcal{X}'$ has the weak Hilbert property over $S$.
\end{proposition}
\begin{proof}
	Let $\mathcal{U}\subseteq \mathcal{X}$ and $\mathcal{U}'\subseteq \mathcal{X}'$ be dense open subschemes such that the morphism $\mathcal{X}' \to \mathcal{X}$ induces an isomorphism $\mathcal{U}' \to \mathcal{U}$. If $k$ is the function field of $S$, then the induced map
    $\mathcal{X}'(S)^{(1)} \cap \mathcal{U}'(k) \to \mathcal{X}(S)^{(1)} \cap \mathcal{U}(k)$ is surjective.
    Indeed, for a near-integral point $s\in \mathcal{X}(S)^{(1)} \cap \mathcal{U}(k)$, there exists a point $s' \in \mathcal{U}'(k)$ lying over $s$, and the valuative criterion of properness shows that $s' \in \mathcal{X}'(S)^{(1)}$.
    Thus, we have $\pi^{-1}\left(\mathcal{X}(S)^{(1)} \cap \mathcal{U}(k)\right) = \mathcal{X'}(S)^{(1)} \cap \mathcal{U'}(k)$.
    Since the intersection of a dense subset and a dense open subset is dense, it is clear that $(\mathcal{X}_k, \mathcal{X}(S)^{(1)})$ is weakly Hilbertian if and only if $(\mathcal{X}_k, \mathcal{X}(S)^{(1)} \cap \mathcal{U}(k))$ is weakly Hilbertian (and analogously for $\mathcal{X}'$).
    Therefore, since $\mathcal{X}'_k \to \mathcal{X}_k$ is a birational morphism of smooth $k$-varieties,
    the claim follows from Lemma~\ref{Lem:HPairsBirational}.
\end{proof}

\subsection{Smooth proper images}

It is shown in \cite[Theorem~3.7]{CDJLZ} that the weak Hilbert property for smooth proper varieties is inherited by images of smooth proper morphisms.
The proof immediately generalizes to the setting of weakly Hilbertian pairs for smooth integral varieties, which readily implies that the weak Hilbert property descends along smooth proper morphisms of regular arithmetic schemes (as regularity implies smoothness of the generic fibres in characteristic $0$).
However, the following more general result given in \cite[Lemma~5.3]{BSFP23BaseChange} allows for a stronger conclusion on arithmetic schemes.

\begin{lemma}[{\cite{BSFP23BaseChange}}]\label{Lem:SmoothImagesBSFP}
	Let $\varphi \colon X \to Y$ be a smooth surjective morphism of normal integral varieties over a field $k$ of characteristic $0$
    and $T \subseteq Y(k)$ a strongly thin subset.
    If the generic fibre of $\varphi$ is geometrically irreducible or if $\varphi$ is proper, then the preimage of $T$ in $X(k)$ is strongly thin.
\end{lemma}

As a direct consequence, we find that the weak Hilbert property descends along proper morphisms of arithmetic schemes that are smooth on the generic fibre.

\begin{proposition}\label{Prop:WHP_descends_smooth_proper}
    Let $S$ be a normal noetherian integral scheme with function field $k$ of characteristic~$0$ and
	let $\varphi\colon \mathcal{X} \to \mathcal{Y}$ be a surjective morphism of arithmetic schemes over $S$ such that
    $\varphi_k \colon \mathcal{X}_k \to \mathcal{Y}_k$ is smooth.
    Assume that $\varphi_k$ is proper or its generic fibre is geometrically irreducible.
    If $\mathcal{X}$ has the weak Hilbert property over $S$, then $\mathcal{Y}$ has the weak Hilbert property over~$S$.
\end{proposition}
\begin{proof}
    Define $X = \mathcal{X}_k, Y = \mathcal{Y}_k$ and  $T = \varphi(\mathcal{X}(S)^{(1)}) \subseteq \mathcal{Y}(S)^{(1)}$.
    Since $\mathcal{X}(S)^{(1)}$ is contained in the preimage of $T$ under the base change $\varphi_k \colon X \to Y$ and not strongly thin in $X$ by the weak Hilbert property of $\mathcal{X}$, Lemma~\ref{Lem:SmoothImagesBSFP} shows that $T$ is not strongly thin in $Y$. In particular, $\mathcal{Y}$ satisfies the weak Hilbert property over $S$.
\end{proof}

\subsection{Finite étale covers and extending the base ring}

The converse to Proposition~\ref{Prop:WHP_descends_smooth_proper} is false even for finite étale morphisms of varieties (cf. \cite[Remark~3.5]{CDJLZ}).
However, it is true for smooth proper varieties up to a finite extension of the base field \cite[Theorem~3.16]{CDJLZ}.
To prove an analogous statement for arithmetic schemes, we show the following analogue of \cite[Lemma~3.14]{CDJLZ}.

\begin{lemma}\label{Lem:FiniteEtaleCovers}
	Let $L/k$ be a finite extension of fields of characteristic $0$, let $X$ be a normal integral variety over $L$, let $Y$ be a normal integral variety over $k$, and let $\pi \colon X \to Y$ be a finite étale morphism of $k$-varieties. Furthermore, let $\Sigma_X \subseteq X(L)$ and $\Sigma_Y \subseteq Y(k)$ be subsets such that $\Sigma_Y \subseteq \pi(\Sigma_X)$.
	If $(Y, \Sigma_Y)$ is a weakly Hilbertian pair over $k$, then $(X, \Sigma_X)$ is a weakly Hilbertian pair over $L$.
\end{lemma}
\begin{proof}
    Since $\Sigma_Y$ is not strongly thin in $Y$,
    it follows from \cite[Lemma~2.1]{BSFP23BaseChange} that $\Sigma_Y$ and thus $\pi(\Sigma_X)$ are not strongly thin in $Y_L$. Thus, we may and do assume that $k=L$.
    Let $(\pi_i \colon Z_i \to X)_{i=1}^n$ be a finite collection of ramified covers and note that each $\pi \circ \pi_i \colon Z_i \to Y$ is a ramified cover.
    Since $(Y, \Sigma_Y)$ is Hilbertian, the set $T := \Sigma_Y \setminus\bigcup_{i=1}^n \pi(\pi_i(Z_i(k)))$ is dense in $Y$.
    Let $y \in T$.
    Since $\Sigma_Y \subseteq \pi(\Sigma_X)$, there exists an element $x \in \Sigma_X$ such that
    $\pi(x) = y$, and clearly $x \notin \pi_i(Z_i(k))$ for any $i$, as otherwise $y \in \pi(\pi_i(Z_i(k)))$.
    Thus, we have
    \[
        T \subseteq \pi \left( \Sigma_X \setminus\bigcup_{i=1}^n \pi_i(Z_i(k)) \right)
    .\]
    Since $T$ is dense in $Y$ and $\pi$ is finite surjective, we conclude that
    $\Sigma_X \setminus\bigcup_{i=1}^n \pi_i(Z_i(k))$ is dense in $X$,
    so $(X, \Sigma_X)$ is weakly Hilbertian, as claimed.
\end{proof}

Applying Lemma~\ref{Lem:FiniteEtaleCovers} to finite étale covers of arithmetic schemes, we obtain the following proposition.

\begin{proposition}\label{Prop:PersistencePotentialWHP}
    Let $S' \to S$ be a finite surjective morphism of normal integral noetherian schemes with function field of characteristic $0$, let $\mathcal{X}$ be an arithmetic scheme over $S'$, let $\mathcal{Y}$ be an arithmetic scheme over $S$,
    and let $\pi \colon \mathcal{X} \to \mathcal{Y}$ be a finite étale morphism of $S$-schemes.
    If $\mathcal{Y}$ has the weak Hilbert property over $S$ and $\mathcal{Y}(S)^{(1)} \subseteq \pi(\mathcal{X}(S')^{(1)})$, then $\mathcal{X}$ has the weak Hilbert property over $S'$.
\end{proposition}
\begin{proof}
    Immediate from Lemma~\ref{Lem:FiniteEtaleCovers}.
\end{proof}

With the same arguments as in the proof of \cite[Theorem~3.16]{CDJLZ}, we are now able to show the analogous statement that the weak Hilbert property ascends along finite étale morphisms of arithmetic schemes, up to a finite étale cover of the base.

\begin{theorem}
    Let $S$ be a regular integral arithmetic scheme and $\varphi \colon \mathcal{X} \to \mathcal{Y}$ a finite étale morphism of arithmetic schemes over $S$.
    If $\mathcal{Y}$ satisfies the weak Hilbert property over $S$
    and $\mathcal{X}_{K(S)}$ is geometrically integral,
    then there exists a finite étale morphism $S'\to S$
    of regular integral arithmetic schemes
    such that $\mathcal{X}_{S'}$ is an arithmetic scheme over $S'$ with the weak Hilbert property over $S'$.
\end{theorem}
\begin{proof}
    By Theorem~\ref{Thm:ChevalleyWeilIntegral}, there exists a finite étale morphism $S' \to S$ of regular integral arithmetic schemes such that
    $\mathcal{Y}(S)^{(1)} \subseteq \varphi(\mathcal{X}(S')^{(1)})$.
    Since $\mathcal{X}_{K(S')}$ is normal integral, $\mathcal{X}_{S'}$ is an arithmetic scheme over $S'$.
    Let $\pi \colon \mathcal{X}_{S'} \to \mathcal{Y}$ be the composition
    $\mathcal{X}_{S'} \to \mathcal{X} \to \mathcal{Y}$. Then $\pi$ is finite étale (hence surjective),
    so $\mathcal{X}_{S'}$ has the weak Hilbert property over $S'$ by Proposition~\ref{Prop:PersistencePotentialWHP}.
\end{proof}

As a further corollary to Proposition~\ref{Prop:PersistencePotentialWHP} and an analogue to \cite[Proposition~3.15]{CDJLZ},
it is easy to see that the weak Hilbert property of an arithmetic scheme $\mathcal{X}$
persists under base change along a finite étale morphism $S' \to S$ by applying Proposition~\ref{Prop:PersistencePotentialWHP} to the morphism $\mathcal{X}_{S'} \to \mathcal{X}$.
However, \cite[Proposition~3.2]{BSFP23BaseChange} already implies the following stronger result.

\begin{proposition}\label{Prop:WHPFlatBaseChange}
    Let $S' \to S$ be a finite type flat morphism of normal integral noetherian schemes with function field of characteristic $0$.
    Let $\mathcal{X}$ be an arithmetic scheme over $S$.
    If $\mathcal{X}$ has the weak Hilbert property over $S$, then $\mathcal{X}_{S'}$ has the weak Hilbert property over $S'$.
\end{proposition}
\begin{proof}
    Note first that $\mathcal{X}_{S'}$ is indeed an arithmetic scheme over $S'$ since $\mathcal{X}_{K(S)}$ is geometrically integral by \cite[Tag~0G69]{stacks-project}.
    Since finite type flat morphisms are open, the image $S_0$ of $S'$ in $S$ is a dense open. Clearly, $\mathcal{X}_{S_0}$ has the weak Hilbert property over $S_0$, since the set of near-integral points can only become larger.
    By \cite[Corollaire~6.1.4]{EGAIV2}, pullback along $S' \to S_0$ preserves codimension of closed subsets, so the inclusion $\mathcal{X}(S_0)^{(1)} \subseteq \mathcal{X}(S')^{(1)}$ holds.
    Let $k$ and $L$ denote the function fields of $S_0$ and $S'$, respectively, and define
    $T = \mathcal{X}(S')^{(1)}$. Since $T \cap \mathcal{X}_k(k) \supseteq \mathcal{X}(S_0)^{(1)}$ is not strongly thin in $\mathcal{X}_k$, it follows from \cite[Proposition~3.2]{BSFP23BaseChange} that $T$ is not strongly thin in $\mathcal{X}_L$, i.e., $\mathcal{X}_{S'}$ has the weak Hilbert property over $S'$.
\end{proof}

\begin{remark}
    Note that Proposition~\ref{Prop:WHPFlatBaseChange} also holds if $S' \to S$ is flat surjective and $K(S')/K(S)$ is small by repeating the proof and replacing \cite[Proposition~3.2]{BSFP23BaseChange} with \cite[Proposition~4.2]{BSFP23BaseChange}.
\end{remark}

\section{Conjecture~\ref{Conj:CZ_arithmetic} for algebraic groups}\label{Sec:AlgGroups}

In order to prove Theorem~\ref{Thm:AlgGroupsSatisfyConjCZ}, we combine results from \cite{CDJLZ} on abelian varieties and \cite{FeiLiu} on linear algebraic groups using a fibration theorem due to Liu \cite[Proposition~2.15]{FeiLiu}.
We briefly recall some terminology used in \emph{loc.\ cit.}

While the usual Hilbert property can be formulated using the integrality of many fibres, this may not be done for the weak Hilbert property for arbitrary ramified covers (cf. Proposition~\ref{Prop:WHP_IntegralFibres}).
This was first noted by Zannier in \cite{ZannierDuke}, leading to the definition of so-called (PB)-covers.
These are defined in \cite{ZannierDuke} for connected commutative algebraic groups over number fields, in \cite{CDJLZ} for abelian varieties over fields of characteristic $0$, and in \cite{FeiLiu} for connected algebraic groups over finitely generated fields of characteristic $0$. We give the definition provided in \cite{FeiLiu}:

\begin{definition}\label{Def:PB}
    Let $k$ be a finitely generated field of characteristic $0$, let $G$ be a connected algebraic group over $k$, and
    let $Y$ be a normal geometrically irreducible variety over $k$.
    A cover $\pi \colon Y \to G$ satisfies the property (PB) if, for every isogeny $G' \to G_{\overline{k}}$ of connected algebraic groups over $\overline{k}$, the fibre product $G' \times_{G_{\overline{k}}} Y_{\overline{k}}$ is irreducible.
\end{definition}

As is shown in \cite[Lemma~2.8]{FeiLiu} (and similarly in \cite[Proposition~2.1]{ZannierDuke} and \cite[Lemma~4.4]{CDJLZ}), a cover is (PB) if and only if it does not admit any non-trivial étale subcover.

In \cite[Definition~2.9]{FeiLiu}, the following property for subgroups of the rational points of an algebraic group is introduced:
\begin{definition}\label{Def:HIT}
    Let $k$ be a finitely generated field of characteristic $0$, let $G$ be a connected algebraic group over $k$, and let $\Omega \subseteq G(k)$ be a (Zariski) dense subgroup.
    The pair $(G, \Omega)$ satisfies the property \emph{(HIT) (over $k$)} if, for every cover $\pi \colon Y \to G$ satisfying (PB) and every finite field extension $k'/k$, the set of $x\in \Omega$ such that $Y_{x,k'}$ is irreducible is dense in $G$.
\end{definition}

The property (HIT) is known for all pairs $(G, \Omega)$, where $G$ is an abelian variety or a connected linear algebraic group and $\Omega \subseteq G(k)$ is dense:

\begin{proposition}\label{Prop:HIT_linear_or_abelian}
    Let $k$ be a finitely generated field of characteristic $0$, let $G$ be a connected linear algebraic group or an abelian variety over $k$, and let $\Omega \subseteq G(k)$ be a dense subgroup. Then $(G, \Omega)$ satisfies (HIT).
\end{proposition}
\begin{proof}
    Let $Y \to G$ be a (PB)-cover and $k'/k$ a finite field extension.
    By definition, the base change $Y_{k'} \to G_{k'}$ is a (PB)-cover.
    If $G$ is an abelian variety (resp. a connected linear algebraic group), then it follows from
    \cite[Theorem~1.4]{CDJLZ} (resp. \cite[Theorem~1.1]{FeiLiu}) that the set of $x \in \Omega$ for which $Y_{k',x} \cong Y_{x,k'}$ is irreducible is dense.
\end{proof}

Liu shows that the property (HIT) may be stated in terms of finite index cosets if $\Omega$ is finitely generated (similarly to \cite{CDJLZ} for abelian varieties):

\begin{lemma}\label{Lem:HIT_one_cover}
    Let $k$ be a finitely generated field of characteristic $0$, let $G$ be a connected algebraic group over $k$, let $\Omega \subseteq G(k)$ be a finitely generated dense subgroup, and let $Y \to G$ be a (PB)-cover.
    If $(G, \Omega)$ satisfies (HIT), then there exists a finite index coset $C$ of $\Omega$ such that, for all $c \in C$, the fibre $Y_c$ is irreducible.
\end{lemma}
\begin{proof}
    A proof of this result is given in \cite[Lemma~5.2]{FeiLiu} under the assumption that $G$ is linear, but it works in our case too.
    Indeed, Liu only uses the linearity of $G$ to reduce to the case that $Y \to G$ is Galois
    by \cite[Lemma~4.2]{FeiLiu}, which holds for arbitrary algebraic groups (as noted directly above \emph{loc.\ cit.}),
    and to appeal to Chebotarev's density theorem \cite[Theorem~9.11]{SerreNXp}, which
    is stated in \cite[Theorem~B.9]{PinkMumfordTate} in the generality that we need.
\end{proof}

To go from irreducible fibres of (PB)-covers to fibres of ramified covers without rational points, we will work over the étale locus of the considered covers, which the following lemma allows us to do.
This is an extension of \cite[Lemma~4.6]{CDJLZ}, where $G$ is assumed to be an abelian variety. The proof of \emph{loc.\ cit.}\ generalizes to our situation; we repeat it for the sake of completeness.

\begin{lemma}\label{Lem:AvoidClosedSubsets}
    Let $k$ be a finitely generated field of characteristic $0$, let $G$ be a connected algebraic group over $k$, let $\Omega \subseteq G(k)$ be a finitely generated dense subgroup, and let $Z \subsetneq G$ be a closed subscheme. Then there exists a finite index coset $C$ of $\Omega$ such that $C \cap Z = \emptyset$.
\end{lemma}
\begin{proof}
    By the density of $\Omega$, there exists a point $P\in \Omega$ with $P \notin Z$.
    By standard spreading out arguments, there exist
    a $\mathbb{Z}$-finitely generated integral domain~$R$
    with fraction field $k$,
    a connected finite type group scheme $\mathcal{G}$ over $R$ with $\mathcal{G}_k = G$ and $\Omega \subseteq \mathcal{G}(R)$,
    a closed subscheme $\mathcal{Z} \subseteq \mathcal{G}$ with $\mathcal{Z}_k = Z$, and a section $\sigma \colon \Spec R \to \mathcal{G}$ which induces $P$ in $G(k)$.
    Since $P \notin Z$, the preimage $\sigma^*\mathcal{Z}$ is a proper closed subset of $\Spec R$. Thus, there exists a closed point $p \in \Spec R$
    such that $P \bmod p = \sigma(\Spec \kappa(p)) \notin \mathcal{Z}$.
    Since $\kappa(p)$ is finite, $\mathcal{G}(\kappa(p))$ is a finite group, and the kernel of the specialization map $\mathcal{G}(R) \to \mathcal{G}(\kappa(p))$ is a finite index subgroup of $\mathcal{G}(R)$.
    Let $\Omega'$ be the intersection of this kernel with $\Omega$, which is of finite index in $\Omega$, and define the finite index coset $C := P\Omega' \subseteq \Omega$.
    Then, for every $c \in C$, we have $c \equiv P \bmod p$ and therefore $c \bmod p \notin \mathcal{Z}$. Since $Z = \mathcal{Z} \times_{\mathcal{G}} G$, we obtain $c \notin Z$ and thus $C \cap Z = \emptyset$, as claimed.
\end{proof}

\subsection{The property (HIT) for many covers}

We fix the following datum throughout this subsection:
\begin{itemize}
    \item a finitely generated field $k$ of characteristic $0$,
    \item a connected algebraic group $G$ over $k$,
    \item a connected linear algebraic group $H$ over $k$,
    \item an abelian variety $A$ over $k$,
    \item an exact sequence
        \[
            0 \longrightarrow H \longrightarrow G \longrightarrow A \longrightarrow 0
        \]
        of algebraic groups over $k$,
    \item and a finitely generated dense subgroup $\Omega \subseteq G(k)$ such that $\Omega \cap H(k)$ is dense in $H$.
\end{itemize}

The following proposition shows that these assumptions imply the property (HIT) for $(G, \Omega)$ and that they are still satisfied if we replace $\Omega$ by a finite index subgroup.

\begin{proposition}\label{Prop:AssumptionsSubsectionFinIndHIT}
    The pair $(G, \Omega)$ satisfies (HIT).
    Moreover, if $\Omega' \subseteq \Omega$ is a finite index subgroup, then $\Omega'$ is dense in $G$ and $\Omega' \cap H(k)$ is dense in $H$, i.e.,
    the assumptions of this subsection are still satisfied after replacing $\Omega$ by $\Omega'$.
    In particular, $(G, \Omega')$ satisfies (HIT).
\end{proposition}
\begin{proof}
    Let $L/k$ be a finite field extension, let $\tilde{G}$ be a connected algebraic group over $L$, and let $\tilde{\Omega} \subseteq \tilde{G}(L)$ be a subgroup with $\tilde{\Omega}$ dense in $\tilde{G}$.
    Assume that there exists an isogeny $\tilde{G} \to H_L$ or an isogeny $\tilde{G} \to A_L$.
    In particular, $\tilde{G}$ is a connected linear algebraic group or an abelian variety over $L$.
    Then, by Proposition~\ref{Prop:HIT_linear_or_abelian}, $(\tilde{G}, \tilde{\Omega})$ satisfies (HIT).
    By \cite[Proposition~2.15]{FeiLiu}, this implies that $(G, \Omega)$ satisfies (HIT).

    Now let $\Omega' \subseteq \Omega$ be a finite index subgroup, so
    there are finitely many elements $\omega_i \in \Omega$ such that $\Omega = \bigcup_i \omega_i \Omega'$.
    By the irreducibility of $G$ and the density of $\Omega$ in $G$, this implies that there is an index $i$ with $\omega_i \Omega'$ dense in $G$. Multiplication with $\omega_i^{-1}$ is an automorphism of $G$, so we conclude that $\Omega' = \omega_i^{-1} \omega_i \Omega'$ is dense in $G$.
    By the same argument applied to $H$ and $\Omega \cap H(k)$, it now suffices to show that
    $\Omega' \cap H(k)$ is of finite index in $\Omega \cap H(k)$.
    To see this, consider the composed morphism
    $\Omega \cap H(k) \to \Omega \to \Omega/\Omega'$. Its kernel is $\Omega' \cap H(k)$, so it induces an injective homomorphism
    $(\Omega \cap H(k)) / (\Omega' \cap H(k)) \hookrightarrow \Omega / \Omega'$.
    Since the right-hand side is finite, so is the left-hand side, which proves the claim.
\end{proof}

Since the weak Hilbert property is defined using many covers (as opposed to only one), in order to relate it to the property (HIT), we first show that (HIT) may be formulated for many covers. We follow the arguments used in \cite[Lemma~4.14]{CDJLZ} on a similar statement for abelian varieties.

\begin{lemma}\label{Lem:HIT_many_covers}
    Let $(\pi_i \colon Y_i \to G)_{i=1}^n$ be a finite family of (PB)-covers.
    Then there exists a finite index coset $C$ of $\Omega$ such that,
    for all $c\in C$ and all $i=1,\ldots,n$, the fibre $Y_{i,c}$ is irreducible.
\end{lemma}
\begin{proof}
    We argue by induction on $n$. If $n=1$, the claim follows from Lemma~\ref{Lem:HIT_one_cover}.
    Now assume that the statement holds for $(\pi_i)_{i=1}^{n-1}$, so there exists a finite index coset $C_{n-1} \subseteq \Omega$ such that, for all $c\in C_{n-1}$ and all $i=1,\ldots,n-1$, the fibre $Y_{i,c}$ is irreducible.
    Choose an element $c_{n-1} \in C_{n-1}$ and a finite index subgroup $\Omega' \subseteq \Omega$ such that $C_{n-1} = c_{n-1} \Omega'$.
    For every $i=1,\ldots,n$, let $\pi_i'\colon Y_i \to G$ be the composition of $\pi_i$ and the translation on $G$ by $c_{n-1}^{-1}$. In particular, for $i=1,\ldots,n-1$ and $c \in \Omega'$, the fibre of $\pi_i'$ over $c$ is irreducible.
    
    Since $(G, \Omega')$ satisfies (HIT) by Proposition~\ref{Prop:AssumptionsSubsectionFinIndHIT}, there exists a finite index coset $C$ of $\Omega'$ such that, for every $c\in C$, the fibre of $\pi_n'$ over $c$ is irreducible.
    It now follows that the finite index coset $c_{n-1}C$ of $\Omega$ satisfies the claim.
\end{proof}

\subsection{From (HIT) to the weak Hilbert property}

We keep the notation and assumptions of the previous subsection.
To show that $(G, \Omega)$ is weakly Hilbertian,
we first consider a single ramified cover of $G$, combining the arguments in \cite[Proposition~4.17]{CDJLZ} with the Chevalley--Weil type lifting in Theorem~\ref{Thm:ChevalleyWeilIntegral}.
To do this, it is crucial to work with near-integral points.
We first show a helpful lemma that a finite étale cover of $G$ is an isogeny of algebraic groups, possibly up to a finite extension of $k$.

\begin{lemma}\label{Lem:EtaleCoverIsogeny}
    Let $\tilde{G}$ be a connected algebraic group over $k$ and let $\lambda \colon X \to \tilde{G}$ be a finite étale morphism of integral $k$-varieties with $X$ geometrically connected over~$k$.
    Then there is a finite field extension $L/k$ such that $X_L$ has the structure of a connected algebraic group over $L$ that makes $\lambda_L \colon X_L \to \tilde{G}_L$ an isogeny.
\end{lemma}
\begin{proof}
    Since $\tilde{G}$ is connected with a $k$-rational point, it is geometrically connected over $k$, so $X$, being finite étale over $\tilde{G}$, is normal and geometrically connected over $k$.
    
    Since $\overline{k}$ is of characteristic $0$, the étale fundamental group $\pi_1^{\mathrm{\acute{e}t}}(\tilde{G}_{\overline{k}})$ is commutative by \cite[Proposition~1.1]{BrionSzamuely}. In particular, every étale cover of $\tilde{G}_{\overline{k}}$ corresponds to a normal subgroup of $\pi_1^{\mathrm{\acute{e}t}}(\tilde{G}_{\overline{k}})$ and is thus Galois.
    Thus, it follows from \cite[Proposition~1.1]{BrionSzamuely} that $X_{\overline{k}}$ carries the structure of an algebraic group over $\overline{k}$ that makes $\lambda_{\overline{k}}$ an isogeny.
    The claim now follows by descending $\lambda_{\overline{k}}$ and the group structure of $\tilde{G}_{\overline{k}}$ to a finite extension of $k$.
\end{proof}

\begin{proposition}\label{Prop:FiniteIndexCosetOneRamCover}
    Let $\pi \colon Y \to G$ be a ramified cover.
    Then there exists a finite index coset $C$ of $\Omega$ such that $Y_c(k) = \emptyset$ for all $c\in C$.
\end{proposition}
\begin{proof}
		Since the claim is trivially true if $Y(k)=\emptyset$, we may and do assume that $Y$ admits a $k$-point and is thus geometrically connected over $k$.
		Moreover, it suffices to prove the claim after replacing $k$ by a finite extension.
		Thus, by Lemma~\ref{Lem:EtaleCoverIsogeny}, up to a finite extension of $k$, we may and do assume that
		there exists a factorization $Y \overset{\mu}{\rightarrow} G' \overset{\lambda}{\rightarrow} G$ of $\pi$
		such that $\lambda$ is an étale subcover of $\pi$ of maximal degree and $G' \to G$ is an isogeny of algebraic groups.
		Note that $\mu$ is a (PB)-cover of $G'$ by \cite[Lemma~2.8]{FeiLiu}.
    By Lemma~\ref{Lem:AvoidClosedSubsets}, up to replacing $\Omega$ by a finite index coset and composing $\pi$ with a translation on $G$ that maps this coset to a subgroup,
    we may and do assume that $\Omega$ is contained in the étale locus of $\pi$.

    By spreading out, there exists a normal $\mathbb{Z}$-finitely generated integral domain $R$ with fraction field $k$, a connected linear finite type group scheme $\mathcal{H}$ over $R$, and a connected finite type group scheme $\mathcal{G}$ over $R$
    such that $\mathcal{H}_k = H$ and $\mathcal{G}_k = G$.
    By \cite[Proposition~4.4]{CorvajaFixedPoints}, there exists a finitely generated subgroup
    $\tilde{\Omega} \subseteq \Omega \cap H(k)$ that is dense in $H$. Enlarging $R$ if necessary, we may and do assume that the finitely many generators of $\Omega$ and $\tilde{\Omega}$
    are contained in the groups $\mathcal{G}(R)$ and $\mathcal{H}(R)$, respectively.
    From this we obtain $\Omega \subseteq \mathcal{G}(R)$ and $\tilde{\Omega} \subseteq \mathcal{H}(R)$.
    By \cite[Lemma~2.13]{FeiLiu}, there exists a commutative diagram
    \[
        \xymatrix{
            0 \ar[r] & H' \ar[d]^h \ar[r] & G' \ar[d]^{\lambda} \ar[r] & A' \ar[d] \ar[r] & 0 \\
            0 \ar[r] & H \ar[r] & G \ar[r] & A \ar[r] & 0 
        }
    \]
    of connected algebraic groups over $k$ such that every row is exact and every column is an isogeny.
    Let $L/k$ be a finite extension such that the kernel $\Gamma := \ker(\lambda_L)$
		of $\lambda_L \colon G'_L \to G_L$    
    consists only of $L$-rational points.
    Enlarging $L$ further if necessary, by Corollary~\ref{Cor:CWT_NearIntegralToRational}, we may and do assume that $\mathcal{G}(R) \subseteq \lambda(G'(L))$ and $\mathcal{H}(R) \subseteq h(H'(L))$.
 
    Since a finite base change of a (PB)-cover is (PB), the cover $\mu_L \colon Y_L \to G'_L$ is (PB).
    For every $\sigma \in \Gamma$, let
    $\mu_{\sigma} = \sigma \circ \mu_L \colon Y_L \to G'_L$ denote the (PB)-cover
    given as the composition of $\mu_L$ and translation on $G'_L$ by $\sigma$.

    Define $\Omega' := \lambda^{-1}(\Omega) \cap G'(L)$. By our choice of $L$,
    the isogeny $\lambda_L$ induces a surjective homomorphism $\Omega' \to \Omega$ with finite kernel.
    Thus, $\Omega'$ is an extension of finitely generated groups and therefore finitely generated.
    Moreover, $\Omega' \cap H'(L)$ is dense in $H'$, as its image in $H_L$ contains $\tilde{\Omega}$ by the choice of $L$.
    Thus, by Lemma~\ref{Lem:HIT_many_covers}, there exists a finite index coset $C'$ of $\Omega'$ such that, for every $c' \in C'$ and every $\sigma\in \Gamma$,
    the fiber $\mu_{\sigma}^{-1}(c')$
    is irreducible.
    Let $C = \lambda(C')$, which by the surjectivity of $\Omega' \to \Omega$ is a finite index coset of $\Omega$.
    
    We claim that $C$ satisfies the claim of the proposition.
    Let $c \in C$; we have to show that $Y_c(k) = \emptyset$.
    If there is no $x \in G'(k)$ with $\lambda(x) = c$, then there is nothing to show, so we assume that there exists an $x \in G'(k)$ such that $\lambda(x) = c$. To prove the claim, it suffices to show that the fiber $Y_x = \mu^{-1}(x)$ has no $k$-rational point.
    By the definition of $C$, there exists a point $c' \in G'_c(L)$ such that $c' \in C'$.
    Since $\Gamma$ acts transitively on the fibers of $G'_L \to G_L$, there exists an element $\sigma \in \Gamma$ such that $\sigma x = c'$ (as $L$-points of $G'_L$).
    In particular, we have $\mu_{\sigma}^{-1}(c') = Y_{x,L}$.
    Since $c' \in C'$, we have that $Y_{x,L}$ is irreducible, and thus $Y_x$ is irreducible.
    By our assumption that $\pi$ is étale over every element of $\Omega$ and $Y_x$ is irreducible, $Y_x$ is the spectrum of a field extension of $\kappa(x) = k$ of degree $\deg \mu$. Since $\pi$ is ramified and $\lambda$ is étale, we have $\deg \mu \geq 2$. This shows that $Y_x(k) = \emptyset$ for every $x \in G'_c(k)$ and thus $Y_c(k) = \emptyset$, as claimed.
\end{proof}

Repeating the induction arguments used in Lemma~\ref{Lem:HIT_many_covers}, we conclude that $(G, \Omega)$ is in fact weakly Hilbertian.

\begin{proposition}\label{Prop:FiniteIndexCosetManyRamCovers}
    Let $(\pi_i \colon Y_i \to G)_{i=1}^n$ be a finite family of ramified covers.
    Then there exists a finite index coset $C$ of $\Omega$ such that, for all $c\in C$ and all $i=1,\ldots,n$, we have $Y_{i,c}(k) = \emptyset$.
    In particular, $(G, \Omega)$ is weakly Hilbertian.
\end{proposition}
\begin{proof}
    The case $n=1$ is true by Proposition~\ref{Prop:FiniteIndexCosetOneRamCover}.
    The proof by induction is now the same as in Lemma~\ref{Lem:HIT_many_covers}.
\end{proof}

\subsection{Proof of Theorem~\ref{Thm:AlgGroupsSatisfyConjCZ}}

To prove Theorem~\ref{Thm:AlgGroupsSatisfyConjCZ}, we now drop the notation fixed in the previous subsections. We first show that we may always reduce to the case of finitely generated subgroups.

\begin{proposition}\label{Prop:FinGenDenseSubgroup}
    Let $k$ be a finitely generated field of characteristic $0$ and let
    \[
        0 \longrightarrow H \longrightarrow G \overset{p}{\longrightarrow} A \longrightarrow 0
    \]
    be an exact sequence of connected algebraic groups over $k$, where $H$ is linear and $A$ is an abelian variety. Let $\Omega \subseteq G(k)$ be a subgroup that is dense in $G$ and such that $\Omega \cap H(k)$ is dense in $H$. Then there exists a finitely generated subgroup $\Omega'$ of $\Omega$ such that $\Omega'$ is dense in $G$ and $\Omega' \cap H(k)$ is dense in $H$.
\end{proposition}
\begin{proof}
    By \cite[Proposition~4.4]{CorvajaFixedPoints}, there exists a finitely generated subgroup
    $\Omega_H = \langle h_1,\ldots,h_m \rangle$ of $\Omega \cap H(k)$ that is dense in $H$.
    Moreover, since $A(k)$ is finitely generated by the Mordell--Weil Theorem and commutative,
    its subgroup $p(\Omega)$ is finitely generated. Let $a_1,\ldots,a_n \in \Omega$ be elements such that $p(\Omega) = \langle p(a_1),\ldots, p(a_n) \rangle$.
    Define $\Omega' = \langle h_1,\ldots,h_m, a_1,\ldots,a_n \rangle \subseteq \Omega$.

    We first claim that $\Omega = \langle a_1,\ldots, a_n \rangle (\Omega \cap H(k))$ (where the inclusion from right to left is obvious).
    Indeed, for any $x \in \Omega$, we find $r_1,\ldots,r_n \in \mathbb{Z}$ such that
    $p(x) = p(a_1)^{r_1} \cdots p(a_n)^{r_n}$.
    Define $y = a_1^{r_1} \cdots a_n^{r_n} \in \langle a_1,\ldots, a_n \rangle \subseteq \Omega$
    and note that $y^{-1}x \in \Omega \cap \ker(p) = \Omega \cap H(k)$, which implies that
    $x \in y(\Omega \cap H(k)) \subseteq \langle a_1,\ldots, a_n \rangle (\Omega \cap H(k))$, as claimed.

    We now show that the closure of $\Omega'$ in $G$ contains $\Omega$, which by the density of $\Omega$ in $G$ concludes the proof.
    Given a subset $V \subseteq G$, let $\cl{V}$ denote the topological closure of $V$ in $G$.
    Then the inclusion
    \[
        \cl{\Omega'}
        \supseteq \cl{\bigcup_{a \in \langle a_1,\ldots, a_n \rangle} a \Omega_H}
        \supseteq \bigcup_{a \in \langle a_1,\ldots, a_n \rangle} \cl{a \Omega_H}
    \]
    of subsets of $G$ holds.
    Since multiplication by an element $a\in G$ on $G$ is a homeomorphism, we have
    \[
        \cl{a \Omega_H}
        = a\cl{\Omega_H}
        = aH
        \supseteq a(\Omega \cap H(k))
    .\]
    This implies that
    \[
        \cl{\Omega'}
        \supseteq \bigcup_{a \in \langle a_1,\ldots, a_n \rangle} a(\Omega \cap H(k)) 
        = \langle a_1,\ldots, a_n \rangle (\Omega \cap H(k)) 
        = \Omega
    ,\]
    as required.
\end{proof}

We can now prove Theorem~\ref{Thm:AlgGroupsSatisfyConjCZ}.

\begin{proof}[Proof of Theorem~\ref{Thm:AlgGroupsSatisfyConjCZ}]
    Let $k$ be the fraction field of $R$ and $G = \mathcal{G}_k$.
    Note that, by \cite[Tag~047N]{stacks-project}, $G$ is smooth over $k$.
    Therefore, by \cite[Theorem~1.1]{ConradCT}, there exists an exact sequence
    \[
        0 \longrightarrow H \longrightarrow G \longrightarrow A \longrightarrow 0
    \]
    of connected algebraic groups over $k$ with $H$ linear and $A$ an abelian variety.
    By \cite[Corollary~18.3]{BorelLinAlgGrp}, $H(k)$ is dense in $H$.
    Thus, by Proposition~\ref{Prop:FinGenDenseSubgroup}, there exists a finitely generated dense subgroup $\Omega$ of $G(k)$ such that $\Omega \cap H(k)$ is dense in $H$.
    Let $\Spec R' \subseteq \Spec R$ be a dense open such that the generators of $\Omega$ are defined over $R'$, i.e., $\Omega\subseteq \mathcal{G}_{R'}(R')$.
    Since the assumptions of the previous sections are satisfied for $G, H, A, \Omega$,
    it follows from Proposition~\ref{Prop:FiniteIndexCosetManyRamCovers} that
    $(G, \Omega)$ is weakly Hilbertian.
    In particular, $(G, \mathcal{G}_{R'}(R'))$ is weakly Hilbertian,
    so $\mathcal{G}_{R'}$ satisfies the weak Hilbert property over~$R'$.
\end{proof}

Theorem~\ref{Thm:AlgGroupsSatisfyConjCZ} and standard spreading out arguments imply the following:

\begin{corollary}\label{Cor:AlgGroupsConjCZ}
    Let $k$ be a finitely generated field of characteristic $0$ and let $G$ be a connected algebraic group over $k$.
    If $G(k)$ is dense in $G$, then $G$ satisfies the weak Hilbert property over $k$.
\end{corollary}

\bibliography{references}{}
\bibliographystyle{alpha}

\end{document}